\documentclass[aos,MSNbibl,nameyear,seceqn,dvips]{arximspdf}
\usepackage[mathscr]{eucal}
\usepackage{mathbh}

\doi{10.1214/12-AOS1054} %kopijuoti is PTS
\volume{40}
\issue{6}
\pubyear{2012}
\firstpage{2877}
\lastpage{2909}

\makeatletter
\newtheorem{theorem}{Theorem}[section]
\newtheorem{proposition}[theorem]{Proposition}
\newtheorem{lemma}[theorem]{Lemma}
\newtheorem{corollary}[theorem]{Corollary}
\newproclaim{remark}{Remark}[section]

\newcommand{\eqref}[1]{(\ref{#1})}

\def\eqdef{\stackrel{\operatorname{def}}{=}}

\newcommand{\cc}[1]{\mathscr{#1}}

\renewcommand{\tilde}[1]{\widetilde{#1}}

\mathchardef\varDelta="0101
\mathchardef\varTheta="0102
\mathchardef\varUpsilon="0107
\mathchardef\varPsi="0109

\renewcommand{\Delta}{\varDelta}
\renewcommand{\Psi}{\varPsi}
\renewcommand{\Theta}{\varTheta}

\def\Var{\operatorname{Var}}
\def\argmax{\arg\max}
\def\argmin{\arg\min}

\def\tr{\operatorname{tr}}
\def\R{\mathbb{R}}
\def\E{\mathbb{E}}
\def\P{\mathbb{P}}
\def\T{\top}
\def\loc{\operatorname{loc}}

\def\fv{\mathbf{f}}
\def\uv{\mathbf{u}}

\def\Yv{\mathbf{Y}}

\def\gammav{\bolds{\gamma}}
\def\varepsilonv{\bolds{\varepsilon}}
\def\xiv{\bolds{\xi}}
\def\zetav{\bolds{\zeta}}

\def\rdl{\epsilon}
\def\rd{\bolds{\rdl}}
\def\rddelta{\delta}
\def\rdomega{\varrho}

\def\rdb{\rd}
\def\rdm{\underline{\rdb}}

\def\Id{\mathbf{I}}
\def\Ind{\mathbh{1}}

\def\muc{\mu_{c}}
\def\gmd{\gm_{0}}

\def\nsize{{n}}

\def\rhor{\omega}

\def\La{\mathbb{L}}
\def\Lab{\La_{\rdb}}
\def\Lam{\La_{\rdm}}

\def\DP{D}
\def\DPc{\DP_{0}}
\def\DPb{\DP_{\rdb}}
\def\DPm{\DP_{\rdm}}

\def\gmi{\mathtt{b}}
\def\gmiid{\mathtt{g}_{1}}
\def\kullbi{\Bbbk}
\def\Thetasi{\Theta_{\loc}}
\def\rri{\mathtt{u}}
\def\rris{\rri_{0}}

\def\IF{\mathbf{f}}
\def\IFc{\IF_{0}}
\def\IFb{\IF_{\rdb}}

\def\xivb{\xiv_{\rdb}}
\def\xivm{\xiv_{\rdm}}

\def\pen{\mathfrak{t}}

\def\ex{\mathrm{e}}
\def\entrl{\mathbb{Q}}
\def\entrlb{\entrl}

\def\kullb{\cc{K}} %{\wp}

\def\gm{\mathtt{g}}
\def\gmc{\gm_{c}}
\def\gmb{\gm}

\def\yyc{\mathtt{y}_{c}}
\def\xx{\mathtt{x}}
\def\xxc{\xx_{c}}

\def\alp{\alpha}

\def\lossp{u}

\def\bern{q}

\def\rups{\rr_{0}}

\def\PDOM{\bolds{\mu}_{0}}

\def\CS{\cc{E}}

\def\nunu{\nu_{0}}

\def\cdimb{\mathfrak{c}_{1}}

\def\rdomega{\varrho}

\def\err{\diamondsuit}
\def\errm{\err_{\rdm}}
\def\errb{\err_{\rdb}}

\def\LCS{\cc{C}}

\def\LL{\cc{L}}

\def\YY{\cc{Y}}

\def\pnnd{\mathfrak{u}}

\def\dimp{p}
\def\dimA{\mathtt{p}_{0}}
\def\BB{\mathbb{B}}
\def\vA{\mathtt{v}}

\def\spread{\Delta}

\def\thetav{\bolds{\theta}}

\def\Ups{\varUpsilon}

\def\Thetas{\Theta_{0}}

\def\fs{f}

\def\UP{\cc{U}}

\def\VP{V}
\def\VPc{\VP_{0}}

\def\VV{H}

\def\VVc{\VV_{0}}

\def\fis{\mathfrak{a}}

\def\Ldrift{M}

\def\mubc{\mu}

\def\expzeta{\mathfrak{s}}

\def\rr{\mathtt{r}}

\def\zz{\mathfrak{z}}
\def\zzQ{\zz_{0}}

\def\LCS{C}

\def\LCS{C}

\def\AssId{\mathcal{I}}

\def\dimA{\mathtt{p}}

\def\dens{f}

\makeatother

\begin{document}
\begin{frontmatter}

\title{Parametric estimation. Finite sample theory}
\runtitle{Parametric estimation. Finite sample theory}

\begin{aug}
\author[A]{\fnms{Vladimir} \snm{Spokoiny}\corref{}\ead[label=e1]{spokoiny@wias-berlin.de}\thanksref{t1}}
\thankstext{t1}{Supported by Predictive Modeling Laboratory, MIPT, RF government grant, ag. 11.G34.31.0073.}
\runauthor{V. Spokoiny}
\affiliation{Weierstrass-Institute, Humboldt University Berlin and
Moscow~Institute~of~Physics~and Technology}
\address[A]{Weierstrass-Institute\\
Mohrenstr. 39, 10117 Berlin\\
Germany\\
and\\
Humboldt University Berlin\\
Germany\\
and\\
Moscow Institute of Physics and Technology\\
\printead{e1}} %adresu isvedimo komanda gale!
\end{aug}

% HISTORY:
\received{\smonth{11} \syear{2011}}
\revised{\smonth{8} \syear{2012}}

% ABSTRACT
%
\begin{abstract}
The paper aims at reconsidering the famous Le Cam LAN theory.
The main features of the approach which make it different from the
classical one are as follows:
(1) the study is nonasymptotic, that is, the sample size is fixed and
does not tend
to infinity;
(2) the parametric assumption is possibly misspecified and the
underlying data
distribution can lie beyond the given parametric family.
These two features enable to bridge the gap between parametric and nonparametric
theory and to build a unified framework for statistical estimation.
The main results include large deviation bounds for the (quasi) maximum
likelihood and the local quadratic bracketing of the log-likelihood process.
The latter yields a number of important corollaries for statistical
inference: concentration, confidence and risk bounds, expansion of the maximum
likelihood estimate, etc.
All these corollaries are stated in a nonclassical way admitting a
model misspecification
and finite samples.
However, the classical asymptotic results including the efficiency
bounds can be easily
derived as corollaries of the obtained nonasymptotic statements.
At the same time, the new bracketing device works well in the
situations with large or
growing parameter dimension in which the classical parametric theory fails.
The general results are illustrated for the i.i.d. setup as well as for
generalized
linear and median estimation.
The results apply for any dimension of the parameter space and provide
a quantitative
lower bound on the sample size yielding the root-n accuracy.
\end{abstract}

% KEYWORDS

\end{frontmatter}

%s1 #&#
\section{Introduction}\label{sec1}

One of the most popular approaches in statistics is based on the parametric
assumption (PA) that the distribution $ \P$ of the observed data $
\Yv$
belongs to a given parametric family
$ (\P_{\thetav}, \thetav\in\Theta\subseteq\R^{\dimp}) $,
where $ \dimp$
stands for the number of parameters.
This assumption allows to reduce the problem of statistical
inference about $ \P$ to recovering the parameter $ \thetav$.
The theory of parameter estimation and inference is nicely developed in
a quite general
setup.
There is a vast literature on this issue.
We only mention the book by \citet{IH1981}, which provides a comprehensive
study of asymptotic properties of maximum likelihood and Bayesian estimators.
The theory is essentially based on two major assumptions:
(1) the underlying data distribution follows the PA;
(2) the sample size or the amount of available information is large
relative to the
number of parameters.

In many practical applications, both assumptions can be very
restrictive and
limit the scope of applicability for the whole approach.
Indeed, the PA is usually only an approximation of real data
distribution and in most
statistical problems it is too restrictive to assume that the PA is exactly
fulfilled.
Many modern statistical problems deal with very complex
high-dimensional data where a huge number of parameters are involved.
In such situations, the applicability of large sample asymptotics is
questionable.
These two issues partially explain why the parametric and nonparametric
theory are almost
isolated from each other.
Relaxing these restrictive assumptions can be viewed as
an important challenge of the modern statistical theory.
The present paper attempts at developing a unified approach which does not
require the restrictive parametric assumptions but still enjoys the
main benefits of the parametric
theory.

The main steps of the approach are similar to the classical local asymptotic
normality (LAN) theory [see, e.g., Chapters 1--3 in the monograph
\citet{IH1981}]:
first one localizes the problem to a neighborhood of the target parameter.
Then one uses a local quadratic expansion of the log-likelihood to
solve the corresponding
estimation problem.
There is, however, one feature of the proposed approach which makes it
essentially
different from the classical scheme.
Namely, the use of the bracketing device instead of classical Taylor expansion
allows to consider much larger local neighborhoods than in the LAN theory.
More specifically, the classical LAN theory effectively requires a
strict localization to a root-n
vicinity of the true point.
At this point, the LAN theory fails in extending to the nonparametric situation.
Our approach works for any local vicinity of the true point.
This opens the door to building a unified theory including most of the classical
parametric and nonparametric results.

Let $ \Yv$ stand for the available data.
Everywhere below we assume that the observed data $ \Yv$ follow the
distribution
$ \P$ on a metric space $ \YY$.
We do not specify any particular structure of $ \Yv$.
In particular, no assumption like independence or weak
dependence of individual observations is imposed.
The basic parametric assumption is that $ \P$ can be approximated by a
parametric distribution $ \P_{\thetav} $ from a given parametric
family $
(\P_{\thetav}, \thetav\in\Theta\subseteq\R^{\dimp}) $.
Our approach allows that the PA can be misspecified, that is, in general,
$ \P\notin ( \P_{\thetav}  ) $.

Let $ L(\Yv,\thetav) $ be the log-likelihood for the considered
parametric model:
$ L(\Yv,\thetav) = \log\frac{d\P_{\thetav}}{d\PDOM}(\Yv) $,
where $ \PDOM$ is
any dominating measure for the family $ (\P_{\thetav}) $.
We focus on the properties of the process $ L(\Yv,\thetav) $ as a
function of the
parameter $ \thetav$.
Therefore, we suppress the argument $ \Yv$ there and
write $ L(\thetav) $ instead of $ L(\Yv,\thetav) $.
One has to keep in mind that $ L(\thetav) $ is random and depends on the
observed data~$ \Yv$.
By $ L(\thetav,\thetav^{*}) \eqdef L(\thetav) - L(\thetav^{*}) $ we denote
the log-likelihood ratio.
The classical likelihood principle suggests to estimate $ \thetav$ by
maximizing the corresponding log-likelihood function $ L(\thetav) $:
%
%e1.1 #&#
\begin{equation}
%[c]
\tilde{\thetav} \eqdef \mathop{\argmax}_{\thetav\in\Theta} L(\thetav) .
\label{tthetamk}
\end{equation}
Our ultimate goal is to study the properties of the \emph
{quasi-maximum likelihood
estimator} (MLE) $ \tilde{\thetav} $. It turns out that such
properties can be
naturally described in terms of the maximum of the process $ L(\thetav
) $ rather than
the point of maximum $ \tilde{\thetav} $. To avoid technical
burdens, it is assumed
that the maximum is attained leading to the identity $ \max_{\thetav
\in\Theta} L(\thetav) =
L(\tilde{\thetav}) $.
However, the point of maximum does not have to be unique. If
there are many such points, we take $ \tilde{\thetav} $ as any of
them. Basically, the
notation $ \tilde{\thetav} $ is used for the identity $ L(\tilde
{\thetav}) =
\sup_{\thetav\in\Theta} L(\thetav) $.

If $ \P\notin ( \P_{\thetav}  ) $, then
the (quasi) MLE $ \tilde{\thetav} $ from \eqref{tthetamk} is
still meaningful and it appears to be an estimator of the value $
\thetav^{*}$
defined by maximizing the expected value of $ L(\thetav) $:
%
%e1.2 #&#
\begin{equation}
%[c]
\thetav^{*} \eqdef \mathop{\argmax}_{\thetav\in\Theta} \E L(\thetav),
\label{thetavsd}
\end{equation}
which is the true value in the parametric situation and can be viewed
as the
parameter of the best parametric fit in the general case.

The results below show that the main properties of the quasi-MLE
$ \tilde{\thetav} $ like concentration or coverage probability can
be described
in terms of the \emph{excess} which is the difference between the
maximum of the process $ L(\thetav) $ and
its value at the ``true'' point $ \thetav^{*}$:
%
%e1.3 #&#
\[
%[c]
L\bigl(\tilde{\thetav},\thetav^{*}\bigr) \eqdef L(\tilde{
\thetav}) - L\bigl(\thetav^{*}\bigr) = \max_{\thetav\in\Theta} L(\thetav) -
L\bigl(\thetav^{*}\bigr). \label{Excessdef}
\]
The established results can be split into two big groups.
A large deviation bound states some concentration properties of the
estimator $ \tilde{\thetav} $.
For specific local sets $ \Thetas(\rr) $ with elliptic shape,
the deviation probability
$ \P ( \tilde{\thetav} \notin\Thetas(\rr)  ) $
is exponentially small in $ \rr$.
This concentration bound allows to restrict the parameter space to a properly
selected vicinity $ \Thetas(\rr) $.
Our main results concern the local properties of the process $
L(\thetav) $ within
$ \Thetas(\rr) $
including a bracketing bound and its corollaries.

The paper is organized as follows.
Section~\ref{Scondgllo} presents the list of conditions which are
systematically
used in the text.
The conditions only concern the properties of the quasi-log-likelihood
process $ L(\thetav) $.
Section~\ref{SLBuLI} appears to be central in the whole approach and
it focuses on local properties of the process $ L(\thetav) $
within $ \Thetas(\rr) $.
The idea is to sandwich the underlying (quasi) log-likelihood process
$ L(\thetav) $ for $ \thetav\in\Thetas(\rr) $ between two
quadratic (in parameter) expressions.
Then the maximum of $ L(\thetav) $ over $ \Thetas(\rr) $ will be
sandwiched as well
by the maxima of the lower and upper processes.
The quadratic structure of these processes helps to compute these
maxima explicitly
yielding the bounds for the value of the original problem.
This approximation result is used to derive a number of corollaries
including the concentration and coverage probability, expansion of the estimator
$ \tilde{\thetav} $, polynomial risk bounds, etc.
In contrary to the classical theory, all the results are nonasymptotic
and do not involve any small values of the form $ o(1) $,
all the terms are specified explicitly.
Also, the results are stated under possible model misspecification.

Section~\ref{Chgparam} accomplishes the local results with the
concentration property
which bounds the probability
that $ \tilde{\thetav} $ deviates from the local set $ \Thetas
(\rr) $.
In the modern statistical literature there are a number of studies considering
maximum likelihood or, more generally, minimum contrast estimators in a general
i.i.d. situation, when the parameter set $ \Theta$ is a subset of some
functional space.
We mention the papers of \citet{vdG1993},
Birg{\'e} and Massart (\citeyear{BiMa1993,BiMa1998}), \citet{Bi2006} and the
references therein.
The established results are based on deep probabilistic facts from empirical
process theory; see, for example,
Talagrand (\citeyear{Ta1996,Ta2001,Ta2005}),
\citet{VW1996} and
\citet{BuLuMa2003}.
The general result presented in Section~2 of the supplement [\citet{supp}] follows the
generic chaining
idea due to \citet{Ta2005}; cf. \citet{Be2006}.
However, we do not assume any specific structure of the model.
In particular, we do not assume independent observations and, thus,
cannot apply
the most developed concentration bounds from the empirical process theory.

Section~\ref{Sexpex} illustrates the applicability of the general results
to the classical case of an i.i.d. sample.
The previously established general results apply under rather mild conditions.
Basically we assume some smoothness of the
log-likelihood process and some minimal number of observations per parameter:
the sample size should be at least of order of the dimensionality $
\dimp$
of the parameter space.
We also consider the examples of generalized linear modeling and of median
regression.

It is important to mention that the nonasymptotic character of our
study yields an
almost complete change of the mathematical tools: the notions of
convergence and
tightness become meaningless, the arguments based on compactness of the
parameter
space do not apply, etc.
Instead we utilize the tools of the empirical process theory based on
the ideas of
concentration of measures and nonasymptotic entropy bounds.
Section~2 of the supplement [\citet{supp}]
presents an exponential bound for a general quadratic form
which is very important for getting the sharp risk bounds for the quasi-MLE.
This bound is an important step in the concentration results for the quasi-MLE.
Section~1 of the supplement [\citet{supp}] explains how the generic chaining and
majorizing measure device
by \citet{Ta2005} refined in \citet{Be2006}
can be used for obtaining a general exponential bound for the log-likelihood
process.

The proposed approach can be useful in many further
research directions including
penalized maximum likelihood and semiparametric estimation [\citet{AASp2012}],
contraction rate and asymptotic normality of the posterior within the
Bayes approach
[\citet{SP2012}]
and local adaptive quantile estimation [\citet{SpWe2012}].

%================================ cond ================================
%s2 #&#
\section{Conditions}
\label{Scondgllo}
Below we collect the list of conditions which are systematically used
in the
text.
It seems to be an advantage of the whole approach that all the results
are stated in a
unified way under the same conditions.
Once checked, one obtains automatically all the established results.
%which are quite general and not very much restrictive.
We do not try to formulate the conditions and the results in the most
general form.
In some cases we sacrifice generality in favor of readability and ease of
presentation.
%If possible we indicate possible extensions of the results under more
%general
%conditions.
It is important to stress that all the conditions only concern the
properties of the quasi-likelihood process $ L(\thetav) $.
Even if the process $ L(\cdot) $ is not a sufficient statistic, the
whole analysis
is entirely based on its geometric structure and probabilistic properties.
The conditions are not restrictive and can be effectively checked in
many particular
situations.
Some examples are given in Section~\ref{Chgexamples} for i.i.d setup,
generalized linear
models and for median regression.

The imposed conditions can be classified into the following groups by their
meaning:
\begin{itemize}
\item smoothness conditions on $ L(\thetav) $ allowing the second
order Taylor expansion;
\item exponential moment conditions;
\item identifiability and regularity conditions.
\end{itemize}
We also distinguish between local and global conditions.
The global conditions concern the global behavior of the process $
L(\thetav) $
while the local conditions focus on its behavior in the vicinity of the central
point $ \thetav^{*}$.
Below we suppose that degree of locality is described by a number $
\rr$.
The local zone corresponds to $ \rr\le\rups$ for a fixed $ \rups
$.
The global conditions concern $ \rr> 0 $.

%s2.1 #&#
\subsection{Local conditions}
\label{Sloccond}
Local conditions describe the properties of $ L(\thetav) $ in a
vicinity of the
central point $ \thetav^{*}$ from \eqref{thetavsd}.

To bound local fluctuations of the process $ L(\thetav) $, we
introduce an
exponential moment condition on the stochastic component $ \zeta
(\thetav) $:
%
%e2.1 #&#
\[
%[c]
\zeta(\thetav) \eqdef L(\thetav) - \E L(\thetav).
\]
Below we suppose that the random function $ \zeta(\thetav) $ is
differentiable in
$ \thetav$ and its gradient $ \nabla\zeta(\thetav) = \partial
\zeta(\thetav) / \partial\thetav\in\R^{\dimp} $ has some
exponential moments.
Our first condition describes the property of the gradient $ \nabla
\zeta(\thetav^{*}) $
at the central point $ \thetav^{*}$.

\begin{longlist}[(${E D_{0}}$)]
\item[(${E D_{0}}$)]
\textit{There exist a positive symmetric matrix $ \VPc^{2} $
and constants $ \gmb> 0 $, $ \nunu\ge1 $ such that
$ \Var \{ \nabla\zeta(\thetav^{*})  \} \le\VPc^{2} $ and
for all
$ |\lambda| \le\gmb$}
%
%e2.2 #&#
\[
%[c]
\label{expzetacloc} \sup_{\gammav\in\R^{\dimp}} %\sup_{\thetav\in\Thetas(\rr)}
\log\E\exp \biggl
\{ \lambda\frac{\gammav^{\T} \nabla\zeta(\thetav^{*})} {
\| \VPc\gammav\|} \biggr\} \le \nunu^{2} \lambda^{2}
/ 2.
\]
\end{longlist}

In a typical situation, the matrix $ \VPc^{2} $ can be defined as
the covariance matrix
of the gradient vector $ \nabla\zetav(\thetav^{*}) $:
$ \VPc^{2} = \Var ( \nabla\zetav(\thetav^{*})  )
= \Var ( \nabla L(\thetav^{*})  ) $.
If $ L(\thetav) $ is the log-likelihood for a correctly specified model,
then $ \thetav^{*}$ is the true parameter value and $ \VPc^{2} $
coincides with
the corresponding Fisher information matrix.
The matrix $ \VPc$ shown in this condition determines the local geometry
in the vicinity of $ \thetav^{*}$. In particular, define the local elliptic
neighborhoods of $ \thetav^{*}$ as
%
%e2.3 #&#
\begin{equation}
%[c]
\Thetas(\rr) \eqdef \bigl\{ \thetav\in\Theta\dvtx \bigl\| \VPc\bigl(
\thetav- \thetav^{*}\bigr)\bigr \| \le\rr\bigr\}. \label{Theta0R}
\end{equation}
The further conditions are restricted to such defined neighborhoods
$ \Thetas(\rr) $.
%In fact, they quantify local smoothness properties of the
%log-likelihood function \( L(\thetav) \).

\begin{longlist}[$({E D_{1}})$]
\item[$({E D_{1}})$]\textit{For each $ \rr\le
\rups$, there exists a
constant $ \rhor(\rr) \le1/2 $ such that it holds for all $
\thetav\in
\Thetas(\rr) $}
%
%e2.4 #&#
\[
\label{expzetacloc1} \sup_{\gammav\in\R^{\dimp}} %\sup_{\thetav\in\Thetas(\rr)}
\log\E\exp \biggl\{ \lambda
\frac{\gammav^{\T} \{ \nabla\zeta(\thetav) - \nabla\zeta
(\thetav^{*}) \}} {
\rhor(\rr) \| \VPc\gammav\|} \biggr\} \le \nunu^{2} \lambda^{2} / 2 ,\qquad |
\lambda| \le\gmb.
\]
Here the constant $ \gmb$ is the same as in $ (E D_{0}) $.
\end{longlist}

The main bracketing result also requires second order smoothness of the
expected log-likelihood $ \E L(\thetav) $.
By definition, $ L(\thetav^{*},\thetav^{*}) \equiv0 $ and $ \nabla
\E
L(\thetav^{*}) = 0 $
because $ \thetav^{*}$ is the extreme point of $ \E L(\thetav) $.
Therefore, $ - \E L(\thetav,\thetav^{*}) $ can be approximated by a quadratic
function of $ \thetav- \thetav^{*}$ in the neighborhood of $ \thetav^{*}$.
The \emph{local identifiability} condition quantifies this quadratic
approximation from above and
from below on the set $ \Thetas(\rr) $ from \eqref{Theta0R}.
%Let the matrix \( \VPc\) be shown in \( (E D_{0}) \).

\begin{longlist}[$({\LL_{0}})$]
\item[$({\LL_{0}})$]
\textit{There is a symmetric strictly positive-definite matrix $ \DPc^{2} $
and for each $ \rr\le\rups$ and a constant $ \rddelta(\rr) \le
1/2 $, such that it holds
on the set $ \Thetas(\rr) = \{ \thetav\dvtx  \| \VPc(\thetav- \thetav^{*}
) \| \le\rr\},
$}
%
%e2.5 #&#
\[
%[c]
\label{LmgfquadEL} \biggl| \frac{- 2 \E L(\thetav,\thetav^{*})}{\| \DPc(\thetav- \thetav^{*}) \|^{2}} - 1 \biggr| \le \rddelta(\rr) .
\]
\end{longlist}

Usually $ \DPc^{2} $ is defined as the negative Hessian of $ \E
L(\thetav^{*}) $:
$ \DPc^{2} = - \nabla^{2} \E L(\thetav^{*}) $.
If $ L(\thetav,\thetav^{*}) $ is the log-likelihood ratio and
$ \P= \P_{\thetav^{*}} $, then
$ - \E L(\thetav,\thetav^{*}) = \break\E_{\thetav^{*}} \log ( d\P_{\thetav^{*}}/ d\P_{\thetav}  )
= \kullb(\P_{\thetav^{*}},\P_{\thetav}) $, the Kullback--Leibler
divergence between
$ \P_{\thetav^{*}} $ and $ \P_{\thetav} $.
Then condition $ (\LL_{0}) $ with $ \DPc= \VPc$ follows from the usual
regularity conditions on the family $ (\P_{\thetav}) $; cf. \citet{IH1981}.
If the log-likelihood process $ L(\thetav) $ is sufficiently smooth
in $ \thetav
$, for example, three times stochastically differentiable,
then the quantities $ \rhor(\rr) $ and $ \rddelta(\rr) $ can be
taken proportional
to the value $ \varrho(\rr) $ defined as
%
%e2.6 #&#
\[
%[c]
\varrho(\rr) \eqdef \max_{\thetav\in\Thetas(\rr)} \bigl\| \thetav-
\thetav^{*}\bigr\| . \label{rhorED1}
\]
In the important special case of an i.i.d. model one can take
$ \rhor(\rr) = \rhor^{*}\rr/ \nsize^{1/2} $ and
$ \rddelta(\rr) = \rddelta^{*}\rr/ \nsize^{1/2} $ for some constants
$ \rhor^{*}, \rddelta^{*}$; see Section~\ref{SqMLEiid}.

The \emph{identifiability condition} relates the matrices
$ \DPc^{2} $ and $ \VPc^{2} $.
\begin{longlist}[$({{\AssId}})$]
\item[$({{\AssId}})$]
There is a constant $ \fis> 0 $ such that
$ \fis^{2} \DPc^{2} \ge\VPc^{2} $.
\end{longlist}

%s2.2 #&#
\subsection{Global conditions}
\label{Sglobalcondition}

The global conditions have to be fulfilled for all $ \thetav$ lying beyond
$ \Thetas(\rups) $.
We only impose one condition on the smoothness of the stochastic
component of the
process $ L(\thetav) $ in term of its gradient and one
identifiability condition in terms
of the expectation $ \E L(\thetav,\thetav^{*}) $.

%The first global condition \( (E) \) assumes some exponential moments
%for the quasi
%log-likelihood \( L(\thetav) \).
%The formulation involves a subset \( \Mubc\) of \( \R_{+} \)
%describing
%all possible exponents in the moment conditions.
%
%For each \( \thetav\in\Theta\), there exists a positive value \(
%such that }
% \E\exp\bigl\{ \mubc L(\thetav,\thetavs) \bigr\}
% < \infty.
%Note that this condition is automatically fulfilled if \( \P= \P_{
%and all the \( \P_{\thetav} \)'s are absolutely
%continuous w.r.t. \( \P_{\thetavs} \) with \( \mubc\le1 \).
%Indeed, \( L(\thetav,\thetavs) = \log\bigl( d\P_{\thetav}/d\P_{
%and \( \log\E_{\thetavs} (d\P_{\thetav}/d\P_{\thetavs}) = 0 \).
%For \( \mubc< 1 \), it holds by the Jensen inequality that
%

The first condition is similar to the local condition $ (E D_{0}) $
and it requires
some exponential moment of the gradient $ \nabla\zeta(\thetav) $
for all
$ \thetav\in\Theta$.
%The global version of \( (E D_{0}) \) assumes that this condition is
%fulfilled for
%any \( \rr\ge\rups\),
However, the constant $ \gm$ may be dependent of the radius
$ \rr= \| \VPc(\thetav- \thetav^{*}) \| $.

\begin{longlist}[$({E\rr})$]
\item[$({E\rr})$]
\textit{For any $ \rr$, there exists a value $
\gm(\rr) > 0 $ such that for all $ \lambda\le\gm(\rr) $}
%
%e2.7 #&#
\[
%[c]
\label{expzetag} \sup_{\thetav\in\Thetas(\rr)} \sup_{\gammav\in\R^{\dimp}} \log\E\exp
\biggl\{ \lambda\frac{\gammav^{\T} \nabla\zeta(\thetav)} {
\| \VPc\gammav\|} \biggr\} \le \nunu^{2}
\lambda^{2} / 2.
\]
\end{longlist}

%Finally we need a global identification property.
%A pointwise identification can be formulated in terms of the Chernoff
%function
%all admissible
%A global concentration bound requires that this result is fulfilled
%uniformly over
%all \( \thetav\).
%We suppose that the deterministic component \( \E L(\thetav,\thetavs)
%log-likelihood is at least comparable with its variance \( \Var L(

The global identification property means that the deterministic component
$ \E L(\thetav,\thetav^{*}) $ of the log-likelihood is competitive
with its
variance $ \Var L(\thetav,\thetav^{*}) $.

\begin{longlist}[$({\cc{L}\rr})$]
\item[$({\cc{L}\rr})$]
\textit{There is a function $ \gmi(\rr) $ such that $ \rr\gmi
(\rr) $
monotonously increases in $ \rr$ and
for each $ \rr\ge\rups$}
% and any \( \thetav\) with \( \| \VPc(\thetav- \thetavs) \| = \rr\)
%
%e2.8 #&#
\[
\inf_{\thetav\dvtx    \| \VPc(\thetav- \thetav^{*}) \| = \rr} \bigl| \E L\bigl(\thetav,\thetav^{*}\bigr) \bigr|
% \frac{- \E L(\thetav,\thetavs)}{\| \VPc(\thetav- \thetavs) \|^{2}}
\ge \gmi(\rr) \rr^{2}. \label{xxentrtt}
\]
%
%for \( b \ge12 \nunu\) and \( b_{1} \ge\).
\end{longlist}

%The function \( \gmi(\rr) \) should not decay too fast with the radius
%Section~\ref{Chgparam} provides sufficient conditions ensuring the
%concentration property
%in terms of the functions \( \gm(\rr) \) and \( \gmi(\rr) \).

%================================ loc ================================
%s3 #&#
\section{Local inference}
\label{Chglocal}
\label{SLBuLI}

The \emph{local asymptotic normality} (LAN) condition since introduced in
\citet{lecam1960} became one of the central notions in the statistical theory.
It postulates a kind of local approximation of the log-likelihood of
the original
model by the log-likelihood of a Gaussian shift experiment.
The LAN property being once checked yields a number of important
corollaries for
statistical inference.
In words, if you can solve a statistical problem for the Gaussian shift
model, the
result can be translated under the LAN condition to the original setup.
We refer to \citet{IH1981} for a nice presentation of the LAN theory including
asymptotic efficiency of MLE and Bayes estimators.
The LAN property was extended to \emph{mixed LAN} or \emph{local asymptotic
quadraticity} (LAQ); see, for example, \citet{LY2000}.
All these notions are very much asymptotic and very much local.
The LAN theory also requires that $ L(\thetav) $ is the correctly specified
log-likelihood.
The strict localization does not allow for considering a growing or infinite
parameter dimension and limits applications of the LAN theory to nonparametric
estimation.

Our approach tries to avoid asymptotic constructions and attempts to
include a
possible model misspecification and a large dimension of the parameter space.
The presentation below shows that such an extension of the LAN theory
can be made
essentially for free:
all the major asymptotic results like Fisher and Cram\'{e}r-Rao
information bounds,
as well as the Wilks phenomenon, can be derived as corollaries of the obtained
nonasymptotic statements simply by letting the sample size to infinity.
At the same time, it applies to a high-dimensional parameter space.

The LAN property states that the considered process $ L(\thetav) $
can be approximated by a quadratic in $ \thetav$ expression
in a vicinity of the central point $ \thetav^{*}$.
This property is usually checked using the second order Taylor expansion.
The main problem arising here is that the error of the approximation
grows too fast with the
local size of the neighborhood.
%This limits the applicability of the methods based on local
%approximation.
Section~\ref{SLqm} presents the nonasymptotic version of the LAN property
in which the local quadratic approximation of $ L(\thetav) $
is replaced by bounding this process from above and from
below by two different quadratic in $ \thetav$ processes.
More precisely, we apply the \emph{bracketing} idea:
the difference $ L(\thetav,\thetav^{*}) = L(\thetav) - L(\thetav^{*}) $
is put between two quadratic processes $ \Lam(\thetav,\thetav^{*}) $ and
$ \Lab(\thetav,\thetav^{*}) $:
%
%e3.1 #&#
\begin{equation}
%[c]
\Lam\bigl(\thetav,\thetav^{*}\bigr) - \errm \le L\bigl(
\thetav,\thetav^{*}\bigr) \le \Lab\bigl(\thetav,\thetav^{*}
\bigr) + \errb, \qquad \thetav\in\Thetas(\rr), \label{LLLerrmb}
\end{equation}
where $ \rdb$ is a numerical parameter, $ \rdm= - \rdb$,
and $ \errm$ and $ \errb$ are stochastic errors which only depend
on the
selected vicinity $ \Thetas(\rr) $.
The upper process $ \Lab(\thetav,\thetav^{*}) $ and the lower process
$ \Lam(\thetav,\thetav^{*}) $ can deviate substantially from each other,
however, the errors $ \errb, \errm$ remain small even if the value
$ \rr$
describing the size of the local neighborhood $ \Thetas(\rr) $ is large.
%and the sandwiching approach applies as long as the local conditions
%are fulfilled.

The sandwiching result \eqref{LLLerrmb} naturally leads to two
important notions:
the {value of the problem} and the {spread}.
It turns out that most of the statements like confidence and concentration
probability rely upon the maximum of $ L(\thetav,\thetav^{*}) $ over
$ \thetav$
which we call \emph{the excess}.
Its expectation will be referred to as \emph{the value of the problem}.
Due to \eqref{LLLerrmb}, the excess can be bounded from above and from below
using the similar quantities $ \max_{\thetav} \Lam(\thetav
,\thetav^{*}) $ and
$ \max_{\thetav} \Lab(\thetav,\thetav^{*}) $ which can be called
the \emph{lower and upper excess}, while their expectations are
\emph{the values of the lower and upper problems}.
Note that
$ \max_{\thetav} \{ \Lab(\thetav,\thetav^{*}) - \Lam(\thetav
,\thetav^{*}) \} $ can be very
large or even infinite.
However, this is not crucial.
What really matters is the difference between the upper and the lower excess.
The \emph{spread} $ \spread_{\rd} $ can be defined as the width of the
interval bounding the excess due to \eqref{LLLerrmb}, that is,
as the sum of the approximation errors and of the difference between
the upper and
the lower excess:
%
%e3.2 #&#
\[
%[c]
\spread_{\rd} \eqdef \errb+ \errm+ \Bigl\{
\max_{\thetav} \Lab\bigl(\thetav,\thetav^{*}\bigr) -
\max_{\thetav} \Lam \bigl(\thetav,\thetav^{*}\bigr) \Bigr\} .
\label{deficdef}
\]
The range of applicability of this approach can be described by the following
mnemonic rule:
``The value of the upper problem is larger in order than the {spread}.''
The further sections explain in detail the meaning and content of this rule.
Section~\ref{SLqm} presents the key bound \eqref{LLLerrmb} and
derives it from the
general results on empirical processes.
Section~\ref{Slocinfr} presents some straightforward corollaries of
the bound
\eqref{LLLerrmb} including the coverage and concentration probabilities,
expansion of the MLE and the risk bounds.
It also %briefly discusses the range of applicability of the local
%approach and
indicates how the classical results on asymptotic
efficiency of the MLE follow from the obtained nonasymptotic bounds.

%s3.1 #&#
\subsection{Local quadratic bracketing}
\label{SLqm}
This section presents
the key result about local quadratic approximation
of the quasi-log-likelihood process given by Theorem~\ref{TapproxLL} below.

Let the radius $ \rr$ of the local neighborhood $ \Thetas(\rr) $
be fixed in
a way that the deviation probability
$ \P ( \tilde{\thetav} \notin\Thetas(\rr)  ) $ is
sufficiently small.
Precise results about the choice of $ \rr$ which ensures this
property are
postponed until Section~\ref{Chgparam}.
In this neighborhood $ \Thetas(\rr) $ we aim at building some
quadratic lower and upper
bounds for the process $ L(\thetav) $.
The first step is the usual decomposition of this process into
deterministic and
stochastic components:
%
%e3.3 #&#
\[
%[c]
L(\thetav) = \E L(\thetav) + \zeta(\thetav), \label{LtvELtvz}
\]
where $ \zeta(\thetav) = L(\thetav) - \E L(\thetav) $.
Condition $ (\LL_{0}) $ allows to approximate the smooth
deterministic function
$ \E L(\thetav) - \E L(\thetav^{*}) $ around the point of maximum $
\thetav^{*}$
by the quadratic form $ - \| \DPc(\thetav- \thetav^{*}) \|^{2}/2 $.
The smoothness properties of the stochastic component $ \zeta(\thetav
) $ given by
conditions $ (E D_{0}) $ and $ (E D_{1}) $ lead to linear approximation
$ \zeta(\thetav) - \zeta(\thetav^{*})
\approx(\thetav- \thetav^{*})^{\T} \nabla\zeta(\thetav^{*}) $.
Putting these two approximations together yields the following
approximation of the
process $ L(\thetav) $ on $ \Thetas(\rr) $:
%
%e3.4 #&#
\begin{equation}
%[c]
L\bigl(\thetav,\thetav^{*}\bigr) \approx \La\bigl(\thetav,
\thetav^{*}\bigr) \eqdef \bigl(\thetav- \thetav^{*}
\bigr)^{\T} \nabla\zeta\bigl(\thetav^{*}\bigr) - \bigl\| \DPc\bigl(
\thetav- \thetav^{*}\bigr) \bigr\|^{2}/2 . \label{LttsttnzDP}
\end{equation}
This expansion is used in most of statistical calculus.
However, it does not suit our purposes because the error of
approximation grows
quadratically with the radius~$ \rr$ and starts to dominate at some critical
value of $ \rr$.
We slightly modify the construction by introducing two different approximating
processes.
They only differ in the deterministic quadratic term which is either
shrunk or
stretched relative to the term $ \| \DPc(\thetav- \thetav^{*}) \|^{2}/2 $
in $ \La(\thetav,\thetav^{*}) $.

Let $ \rddelta, \rdomega$ be nonnegative constants.
Introduce for a vector $ \rd= (\rddelta,\rdomega) $ the following notation:
%
%e3.5 #&#
%e3.6 #&#
\begin{eqnarray}
\label{bLquadloc} \Lab\bigl(\thetav,\thetav^{*}\bigr)& \eqdef &\bigl(
\thetav- \thetav^{*}\bigr)^{\T} \nabla L\bigl(
\thetav^{*}\bigr) - \bigl\| \DPb\bigl(\thetav- \thetav^{*}\bigr)
\bigr\|^{2}/2
\nonumber
\\[-8pt]
\\[-8pt]
\nonumber
&=& \xivb^{\T} \DPb\bigl(\thetav- \thetav^{*}\bigr) -\bigl \| \DPb
\bigl(\thetav- \thetav^{*}\bigr) \bigr\|^{2}/2 ,
\end{eqnarray}
where $ \nabla L(\thetav^{*}) = \nabla\zeta(\thetav^{*}) $ by $
\nabla
\E L(\thetav^{*}) = 0 $
and
%
%e3.7 #&#
\[
\DPb^{2} = \DPc^{2} (1 - \rddelta) - \rdomega
\VPc^{2},\qquad  \xivb \eqdef \DPb^{-1} \nabla L\bigl(
\thetav^{*}\bigr) . \label{xivalpsn}
\]
Here we implicitly assume that with the proposed choice of the constants
$ \rddelta$ and $ \rdomega$,
the matrix $ \DPb^{2} $ is nonnegative: $ \DPb^{2} \ge0 $.
The representation \eqref{bLquadloc} indicates that the process
$ \Lab(\thetav,\thetav^{*}) $ has the geometric structure of
log-likelihood of a
linear Gaussian model.
We do not require that the vector $ \xivb$ is Gaussian and, hence,
it is not the
Gaussian log-likelihood.
However, the geometric structure of this process appears to be more
important than
its distributional properties.

One can see that if $ \rddelta,\rdomega$ are positive, the
quadratic drift
component of the process $ \Lab(\thetav,\thetav^{*}) $ is shrunk
relative to
$ \La(\thetav,\thetav^{*}) $ in \eqref{LttsttnzDP} for
$ \rdb$ positive and it is stretched if $ \rddelta,\rdomega$ are
negative.
Now, given $ \rr$, fix some $ \rddelta\ge\rddelta(\rr) $ and
$ \rdomega\ge3\nunu\rhor(\rr) $ with the value $ \rddelta(\rr
) $ from
condition $ (\LL_{0}) $ and $ \rhor(\rr) $ from condition $ (E
D_{1}) $.
Finally set $ \rdm= - \rdb$, so that
$
\DPm^{2}
=
\DPc^{2} (1 + \rddelta) + \rdomega\VPc^{2}.
$

%th3.1 #&#
\begin{theorem}
\label{TapproxLL}
Assume $ (E D_{1}) $ and $ (\LL_{0}) $.
Let for some $ \rr$ the values
$ \rdomega\ge3 \nunu  \rhor(\rr) $ and $ \rddelta\ge\rddelta
(\rr) $ be such that
$ \DPc^{2} (1 - \rddelta) - \rdomega\VPc^{2} \ge0 $.
%Set \( \rdb= (\rddelta,\rdomega) \), \( \rdm= - \rdb= (-\rddelta,-
Then
%
%e3.8 #&#
\begin{equation}
%[c]
\qquad
\Lam\bigl(\thetav,\thetav^{*}\bigr) - \errm(\rr) \le L
\bigl(\thetav,\thetav^{*}\bigr) \le \Lab\bigl(\thetav,
\thetav^{*}\bigr) + \errb(\rr), \qquad \thetav\in\Thetas(\rr), \label{LttbLtt}
\end{equation}
with $ \Lab(\thetav,\thetav^{*}), \Lam(\thetav,\thetav^{*}) $
defined by
\eqref{bLquadloc}.
The error terms $ \errb(\rr) $ and $ \errm(\rr) $ satisfy the bound
\eqref{errbzrr} from Proposition~\ref{Tdefbound}.
\end{theorem}

The proof of this theorem is given in Proposition~\ref{Tdefbound}.

%re1 #&#
\begin{remark}
\label{RDelta1loc}
%If \( \xx\) is not too big, then the value \( \zzQ(\xx,\entrlb) \) is
%close to
%The bounds \eqref{errbzrr} tells us that
This bracketing bound \eqref{LttbLtt} describes some properties of the
log-likelihood
process and the estimator $ \tilde{\thetav} $ is not shown there.
However, it directly implies most of our inference results.
We therefore formulate \eqref{LttbLtt} as a separate statement.
Section~\ref{Sdeficincy} below presents some exponential bounds on
the error terms $ \errb(\rr) $ and $ \errm(\rr) $.
The main message is that under rather broad conditions, these errors
are small
and have only minor impact on the inference for the quasi-MLE~$ \tilde
{\thetav} $.
%In particular, under some standard regularity conditions they are of
%order
\end{remark}

%s3.2 #&#
\subsection{Local inference} \label{Slocinfr}
This section presents a list of corollaries
from the basic approximation bounds of Theorem~\ref{TapproxLL}. The
idea is to replace
the original problem by a similar one for the approximating upper and
lower models. It
is important to stress once again that all the corollaries only rely on the\vadjust{\goodbreak}
\emph{bracketing result} \eqref{LttbLtt} and the \emph{geometric
structure} of the
processes $ \Lab$ and $ \Lam$. Define the \emph{spread} $
\spread_{\rd}(\rr) $
by
%
%e3.9 #&#
\begin{equation}
\spread_{\rd}(\rr) \eqdef \errb(\rr) + \errm(\rr) + \bigl( \| \xivb
\|^{2} - \| \xivm\|^{2} \bigr)/2 . \label{deficiencypa}
\end{equation}
Here $ \xivb= \DPb^{-1} \nabla L(\thetav^{*}) $ and
$ \xivm= \DPm^{-1} \nabla L(\thetav^{*}) $.
The quantity $ \spread_{\rd}(\rr) $ appears to be the price
induced by our
bracketing device.
Section~\ref{Sdeficincy} below presents some probabilistic bounds on
the spread
showing that it is small relative to the other terms.
All our corollaries below
are stated under conditions of Theorem~\ref{TapproxLL}
and implicitly assume that the spread can be nearly ignored.

%s3.2.1 #&#
\subsubsection{Local coverage probability}
Our first result describes the probability of covering $ \thetav^{*}$
by the random set
%
%e3.10 #&#
\begin{equation}
%[c]
\CS(\zz) = \bigl\{ \thetav\dvtx 2 L(\tilde{\thetav},\thetav) \le\zz
\bigr\} . \label{CSzzpa}
\end{equation}

%co3.2 #&#
\begin{corollary}
\label{Tconfloc}
For any $ \zz> 0 $
%it holds with \( \xivb= \DPb^{-1} \nabla L(\thetavs) \)
%
%e3.11 #&#
\begin{equation}
\P \bigl\{ \CS(\zz) \not\ni\thetav^{*}, \tilde{\thetav} \in \Thetas(
\rr) \bigr\} % &=&
% \P\bigl\{ 2 L(\tilde{\thetav},\thetavs) > \zz,   \tilde{\thetav} \in
% \bigr\}
% \
\le \P
\bigl\{ \| \xivb\|^{2} \ge \zz- \errb(\rr) \bigr\} . \label{coverexploc}
\end{equation}
\end{corollary}

\begin{pf}
The bound \eqref{coverexploc} follows from the upper bound of
Theorem~\ref{TapproxLL}
and the statement \eqref{supLat} of Lemma~\ref{LLbreveloc} below.
\end{pf}

Below [see \eqref{PxivbzzBB}] we also present an exponential bound
which helps to answer a
very important question about a proper choice of the critical value $
\zz$ ensuring
a prescribed covering probability.

%It obviously holds on the set \( \{ \tilde{\thetav} \in\Thetas(\rr)
% \{ \tilde{\thetav} \notin\CAb(z) \}
% & = &
% \Bigl\{
% \sup_{\thetav\notin\CAb(z)} L(\thetav,\thetavs) =
% \sup_{\thetav} L(\thetav,\thetavs)
% \Bigr\}
% \\
% & \subseteq&
% \Bigl\{
% \sup_{\thetav\notin\CAb(z)} \Lab(\thetav,\thetavs) + \errb(\rr) \ge
% \sup_{\thetav} \Lam(\thetav,\thetavs) - \errm(\rr)
% \Bigr\} .
%As \( \Lab(\thetav,\thetavs) \) is a
%quadratic function of \( \DPb(\thetav- \thetavs) \); cf.
%its maximum on the complement \( \CAb^{c}(z) \) of the set
%with \( \gamma= z / \| \xivb\| \).
%This implies for all \( \thetav\notin\CAb(z) \)
% \Lab(\thetav,\thetavs)
% \le
% \gamma\| \xivb\|^{2}
% - \gamma^{2} \| \xivb\|^{2} / 2
% =
% z \| \xivb\| - z^{2} / 2 .
%By Lemma~\ref{LLbreveloc}
%Therefore,
% \{ \tilde{\thetav} \notin\CAb(z) \}
% & \subseteq&
% \bigl\{
% z \| \xivb\| - z^{2} / 2
% \ge
% \| \xivm\|^{2}/2 - \errb(\rr) - \errm(\rr)
% \bigr\}
% \\
% & = &
% \bigl\{ z^{2}/2 - z \| \xivb\| + \| \xivb\|^{2} /2
% \le\spread_{\rd}(\rr)
%% \errb(\rr) + \errm(\rr) + \alpha_{\rd}
% \bigr\}
%and \eqref{PlLdrr} follows.

%s3.2.2 #&#
\subsubsection{Local expansion, Wilks theorem and local concentration}
\label{Slocconc}

Now we show how the bound \eqref{LttbLtt} can be used for obtaining a local
expansion of the quasi-MLE~$ \tilde{\thetav} $.
All our results will be conditioned to the random set
$ \LCS_{\rd}(\rr) $ defined as
%
%e3.12 #&#
\begin{equation}
%[c]
\LCS_{\rd}(\rr) \eqdef \bigl\{ \tilde{\thetav} \in\Thetas(
\rr), \bigl\| \VPc\DPm^{-1} \xivm\bigr\| \le\rr \bigr\} . \label{LCSrdrr}
\end{equation}
The second inequality in the definition of $ \LCS_{\rd}(\rr) $ is related
to the solution of the upper and lower problems (cf. Lemma~\ref{LLbreveloc}):
$ \| \VPc\DPm^{-1} \xivm\| \le\rr$ means
$ \tilde{\thetav}_{\rdm} \notin\Thetas(\rr) $, where
$ \tilde{\thetav}_{\rdm} = \argmin_{\thetav} \Lam(\thetav
,\thetav^{*}) $.

Below in Section~\ref{Sdeficincy} we present some upper bounds on the
value $ \rr$
ensuring a dominating probability of this random set.
The first result can be viewed as a finite sample version of the famous
Wilks theorem.

%The basic idea is to plug \( \tilde{\thetav} \) in place of \( \thetav
%definition of \( \errb(\rr) \).

%co3.3 #&#
\begin{corollary}
\label{TquadWilks}
On the random set $ \LCS_{\rd}(\rr) $ from \eqref{LCSrdrr}, it holds
%
%e3.13 #&#
\begin{equation}
\| \xivm\|^{2}/2 - \errm(\rr) \le L\bigl(\tilde{\thetav},
\thetav^{*}\bigr) \le \| \xivb\|^{2}/2 + \errb(\rr) .
\label{VPctttnL}
\end{equation}
\end{corollary}

The next result is an extension of another prominent asymptotic result,
namely, the
Fisher expansion of the MLE.

%co3.4 #&#
\begin{corollary}
\label{Tquadexploc}
On the random set $ \LCS_{\rd}(\rr) $ from \eqref{LCSrdrr}, it holds
%
%e3.14 #&#
\begin{equation}
%[c]
\bigl\| \DPb \bigl( \tilde{\thetav} - \thetav^{*} \bigr) -
\xivb%\DPc^{-1} \nabla L(\thetavs)
\bigr\|^{2} \le 2 \spread_{\rd}(\rr) .
\label{DPbttt}
\end{equation}
%
%with \( \beta_{\rd}(\rr) \) from \eqref{betaRdlp}.
\end{corollary}

The proof of Corollaries~\ref{TquadWilks} and~\ref{Tquadexploc}
relies on
the solution of the upper and lower problems and it is given below at
the end of this section.

Now we describe\vspace*{1pt} concentration properties of $ \tilde{\thetav} $
assuming that $ \tilde{\thetav} $ is restricted to $ \Thetas(\rr
) $.
More precisely, we bound the probability that
$ \| \DPb(\tilde{\thetav} - \thetav^{*}) \| > z $ for a given $ z >
0 $.
%
%co3.5 #&#
\begin{corollary}
\label{Tdevprob}
For any $ z > 0 $, it holds %with \( \alp_{\rd} \) from
%on the random set \( \LCS_{\rd}(\rr) \)
%
%e3.15 #&#
\[
\label{PlLdrr} \P \bigl\{ \bigl\| \DPb\bigl(\tilde{\thetav} - \thetav^{*}
\bigr)\bigr \| > z, \LCS_{\rd}(\rr) \bigr\} \le \P \bigl\{ \| \xivb\| > z
- \sqrt{2\spread_{\rd}(\rr)} \bigr\} . % \\
% & \le&
% \P\Bigl\{
% \bigl(1 + \sqrt{\alp_{\rd}}\bigr) \| \xivb\| > z
% - \sqrt{2 \errb(\rr) + 2 \errm(\rr)}
% \Bigr\} .
\label{PlLdrra}
\]
%where \( \errbm(\rr) \eqdef\max\bigl\{ \errb(\rr),\errm(\rr) \bigr\}
\end{corollary}

%The result can be slightly refined using the quantity
% \alp_{\rd}
% & \eqdef&
% \| \Id_{\dimp} - \DPb\DPm^{-2} \DPb\|_{\infty}
% =
% \lambda_{\max}\bigl( \Id_{\dimp} - \DPb\DPm^{-2} \DPb\bigr) .
%In words, \( \alp_{\rd} \) is the maximum eigenvalue of the matrix
%
%For any \( z > 0 \), it holds %with \( \alp_{\rd} \) from
%%on the random set \( \LCS_{\rd}(\rr) \)
% \P\bigl\{
% \| \DPb(\tilde{\thetav} - \thetavs) \| > z,
%   \LCS_{\rd}(\rr)
% \bigr\}
% & \le&
% \P\bigl\{
% \bigl(1 + \sqrt{\alp_{\rd}}\bigr) \| \xivb\| > z
% - 2 \sqrt{\errbm(\rr)}
% \bigr\} ,
%where \( \errbm(\rr) \eqdef\max\bigl\{ \errb(\rr),\errm(\rr) \bigr\}
%%The first inequality follows directly from the local expansion
%%Further,
%Indeed, the bound \( \| \xivb\|^{2} - \| \xivm\|^{2} \le\alp_{\rd}
%implies \eqref{PlLdrram} due to
% \sqrt{2 \spread_{\rd}(\rr)}
% & \le&
% 2 \sqrt{\errbm(\rr)} + \sqrt{\alp_{\rd}} \| \xivb\| .
%
An interesting and important question is for which $ \zz$ in \eqref{CSzzpa}
the coverage probability of the event $  \{ \CS(\zz) \ni
\thetav^{*} \} $
or for which $ z $ the concentration probability of the event
$ \{ \| \DPb(\tilde{\thetav} - \thetav^{*}) \| \le z \} $ becomes
close to one.
It will be addressed in Section~\ref{Sdeficincy}.

%s3.2.3 #&#
\subsubsection{A local risk bound}
Below we also bound the moments of the excess $ L(\tilde{\thetav
},\thetav^{*}) $
and of the normalized loss $ \DPb(\tilde{\thetav} - \thetav^{*}) $
when $ \tilde{\thetav} $ is restricted to $ \Thetas(\rr) $.
The result follows directly from Corollaries~\ref{TquadWilks} and \ref
{Tquadexploc}.
%Similarly one can bound the
%localized risk \( \E\closs(\tilde{\thetav},\thetavs) \) for a class
%of loss
%functions \( \closs\).

%co3.6 #&#
\begin{corollary}
\label{TLttr}
For $ \lossp> 0 $
%
%e3.16 #&#
\[
\label{LosspT0R} \E \bigl\{ L^{\lossp}\bigl(\tilde{\thetav},
\thetav^{*}\bigr) \Ind \bigl(\tilde{\thetav} \in\Thetas(\rr) \bigr)
\bigr\} \le \E \bigl[ \bigl\{ \| \xivb\|^{2}/2 + \errb(\rr) \bigr
\}^{\lossp} \bigr] . \label{ELlossp}
\]
Moreover, it holds
%
%e3.17 #&#
\[
\E \bigl\{ \bigl\| \DPb\bigl(\tilde{\thetav} - \thetav^{*}\bigr)
\bigr\|^{\lossp} \Ind \bigl( \LCS_{\rd}(\rr) \bigr) \bigr\} \le \E \bigl[
\bigl\{ \| \xivb\| + \sqrt{2 \spread_{\rd}(\rr)} \bigr\}^{\lossp}
\bigr] . \label{Lrttloc}
\]
\end{corollary}

%s3.2.4 #&#
\subsubsection{Comparing with the asymptotic theory}
This section briefly discusses the relation between the established
nonasymptotic
bounds and the classical asymptotic results in parametric estimation.
This comparison is not straightforward because
the asymptotic theory involves the sample size or noise level as the asymptotic
parameter, while our setup is very general and works even for a
``single'' observation.
Here we simply treat $ \rd= (\rddelta,\rdomega) $ as a small parameter.
This is well justified by the i.i.d. case with $ \nsize$ observations,
where it holds $ \rddelta= \rddelta(\rr) \asymp\sqrt{\rr/\nsize
} $ and
similarly for $ \rdomega$; see Section~\ref{Chgexamples} for more details.
The bounds below in Section~\ref{Sdeficincy} show that the spread
$ \spread_{\rd}(\rr) $ from \eqref{deficiencypa} is small and can be
ignored in the asymptotic calculations.
The results of Corollary~\ref{Tconfloc} through~\ref{TLttr} represent
the desired
bounds in terms of deviation bounds for the quadratic form $ \| \xivb
\|^{2} $.

For better understanding the essence of the presented results, consider first
the ``true'' parametric model with the correctly specified log-likelihood
$ L(\thetav) $.
Then $ \DPc^{2} = \VPc^{2} $ is the total Fisher information matrix.
In the i.i.d. case it becomes $ \nsize\IFc$ where\vadjust{\goodbreak} $ \IFc$ is the
usual Fisher
information matrix of the considered parametric family at $ \thetav^{*}$.
In particular, $ \Var \{ \nabla L(\thetav^{*})  \} = \nsize
\IFc$.
So, if
$ \DPb$ is close to $ \DPc$, then $ \xivb$ can be treated as
the normalized score.
Under usual assumptions, $ \xiv\eqdef\DPc^{-1} \nabla L(\thetav^{*})
$ is the asymptotically standard
normal $ \dimp$-vector.
The same applies to $ \xivb$.
Now one can observe that
Corollaries~\ref{Tconfloc} through~\ref{TLttr}
directly imply most of the classical asymptotic statements.
In particular,
Corollary~\ref{TquadWilks} shows that the twice excess
$ 2 L(\tilde{\thetav},\thetav^{*}) $ is nearly $ \| \xivb\|^{2} $
and thus
nearly $ \chi^{2}_{\dimp} $ (Wilks' theorem).
Corollary~\ref{Tquadexploc} yields the expansion
$ \DPb ( \tilde{\thetav} - \thetav^{*} ) \approx\xivb$
(the Fisher expansion)
and, hence, $ \DPb ( \tilde{\thetav} - \thetav^{*} ) $ is
asymptotically
standard normal.
Asymptotic variance of $ \DPb ( \tilde{\thetav} - \thetav^{*}
) $ is nearly
one, so $ \tilde{\thetav} $ achieves the Cram\'er--Rao efficiency bound
in the asymptotic setup.

\subsection{Spread}
\label{Sdeficincy}
This section presents some bounds on the \emph{spread} $ \spread_{\rd}(\rr) $
from~\eqref{deficiencypa}.
This quantity is random but it can be easily evaluated under the
conditions made.
We present two different results: one bounds the errors
$ \errb(\rr) , \errm(\rr) $, while the other presents a deviation bound
on quadratic forms like $ \| \xivb\|^{2} $.
The results are stated under conditions $ (E D_{0}) $ and $ (E
D_{1}) $
in a nonasymptotic way, so the formulation is quite technical.
An informal discussion at the end of this section explains the typical
behavior of
the spread.
The first result accomplishes the bracketing bound \eqref{LttbLtt}.

%pr3.7 #&#
\begin{proposition}
\label{Tdefbound}
Assume $ (E D_{1}) $.
The error $ \errb(\rr) $ in \eqref{LttbLtt} fulfills
%
%e3.18 #&#
\begin{equation}
%[c]
\P \bigl\{ \rdomega^{-1} \errb(\rr) \ge \zzQ(\xx,\entrlb)
\bigr\} \le \exp ( - \xx ) \label{errbzrr}
\end{equation}
with $ \zzQ(\xx,\entrlb) $ given for $ \gmd= \gmb\nunu\ge3 $ by
%
%e3.19 #&#
\[
%[c]
\zzQ(\xx,\entrlb) \eqdef \cases{ ( 1 + \sqrt{\xx+ \entrlb}
)^{2}, &\quad $\mbox{if } 1 + \sqrt{\xx+ \entrlb} \le\gmd,$ \vspace*{2pt}
\cr
1 + \bigl\{ 2\gmd^{-1} (\xx+ \entrlb) + \gmd \bigr\}^{2}, &\quad  $
\mbox{otherwise} ,$ } \label{PUPxxltt}
\]
where $ \entrlb= \cdimb\dimp$ with $ \cdimb= 2 $ for $ \dimp
\ge2 $
and $ \cdimb= 2.7 $ for $ \dimp= 1 $.
Similarly for $ \errm(\rr) $.
\end{proposition}

%re2 #&#
\begin{remark}
\label{Rgmberrb}
The bound \eqref{errbzrr} essentially depends on the value $ \gmb$ from
condition $ (E D_{1}) $.
The result requires that $ \gmb\nunu\ge3 $.
However, this constant can usually be taken of order $ \nsize^{1/2} $
; see
Section~\ref{Chgexamples} for examples.
If $ \gmb^{2} $ is larger in order than $ \dimp+ \xx$, then
$ \zzQ(\xx,\entrlb) \approx\cdimb\dimp+ \xx$.
\end{remark}

\begin{pf}
Consider for fixed $ \rr$ and $ \rdb= (\rddelta,\rdomega) $ the quantity
%
%e3.20 #&#
\[
\label{Delta1loc} \errb(\rr) \eqdef \sup_{\thetav\in\Thetas(\rr)} \biggl\{ L\bigl(\thetav,
\thetav^{*}\bigr) - \E L \bigl(\thetav,\thetav^{*}\bigr) -
\bigl(\thetav- \thetav^{*}\bigr)^{\T} \nabla L\bigl(
\thetav^{*}\bigr) - \frac{\rdomega}{2} \bigl\| \VPc\bigl(\thetav-
\thetav^{*}\bigr) \bigr\|^{2} \biggr\} .
\]
As $ \rddelta\ge\rddelta(\rr) $, it holds
$ - \E L(\thetav,\thetav^{*}) \ge(1 - \rddelta) \DPc^{2} $ and
$ L(\thetav,\thetav^{*}) - \Lab(\thetav,\thetav^{*}) \le\errb
(\rr) $.
Moreover, in view of $ \nabla\E L(\thetav^{*}) = 0 $, the definition of
$ \errb(\rr) $ can be rewritten as
%
%e3.21 #&#
\[
\label{Delta1loc} \errb(\rr) \eqdef \sup_{\thetav\in\Thetas(\rr)} \biggl\{ \zeta\bigl(
\thetav,\thetav^{*}\bigr) - \bigl(\thetav- \thetav^{*}
\bigr)^{\T} \nabla \zeta \bigl(\thetav^{*}\bigr) -
\frac{\rdomega}{2} \bigl\| \VPc\bigl(\thetav- \thetav^{*}\bigr)
\bigr\|^{2} \biggr\} .\vadjust{\goodbreak}
\]
%
%Similarly define
% \errm(\rr)
% & \eqdef&
% \sup_{\thetav\in\Thetas(\rr)}
% \bigl\{
% L(\thetavs,\thetav) - \E L (\thetavs,\thetav)
% - (\thetavs- \thetav)^{\T} \nabla L(\thetavs)
% - \frac{\rdomega}{2} \| \VPc(\thetav- \thetavs) \|^{2}
% \bigr\}
% \\
% &=&
% \sup_{\thetav\in\Thetas(\rr)}
% \bigl\{
% \zeta(\thetavs,\thetav) - (\thetavs- \thetav)^{\T} \nabla\zeta(
% - \frac{\rdomega}{2} \| \VPc(\thetav- \thetavs) \|^{2}
% \bigr\} .
Now the claim of the theorem can be easily reduced to an exponential
bound for
the quantity $ \errb(\rr) $.
We apply Theorem~2.11 of the supplement [\citet{supp}] %\ref{Tsmoothpenlc}
to the process
%
%e3.22 #&#
\[
%[c]
\UP\bigl(\thetav,\thetav^{*}\bigr) = \frac{1}{\rhor(\rr)}
\bigl\{ \zeta\bigl(\thetav,\thetav^{*}\bigr) - \bigl(\thetav-
\thetav^{*}\bigr)^{\T} \nabla \zeta \bigl(\thetav^{*}
\bigr) \bigr\},\qquad \thetav\in\Thetas(\rr), \label{UPloce}
\]
and $ \VVc= \VPc$.
Condition $ (\CS  D) $ follows from
$ (E D_{1}) $ with the same $ \nunu$ and $ \gmb$ in view of
$ \nabla\UP(\thetav,\thetav^{*})
=  \{ \nabla\zeta(\thetav) - \nabla\zeta(\thetav^{*})
\}
/ \rhor(\rr) $.
So, the conditions of Theorem~2.11 in the supplement [\citet{supp}] %\ref{Tsmoothpenlc}
are fulfilled, yielding
\eqref{errbzrr} in view of $ \rdomega\ge3 \nunu  \rhor(\rr) $.
\end{pf}

Due to the main bracketing result, the local excess
$ \sup_{\thetav\in\Thetas(\rr)} L(\thetav,\thetav^{*}) $ can be
put between
similar quantities for
the upper and lower approximating processes up to the error terms $
\errb(\rr), \errm(\rr) $.
The random quantity $ \sup_{\thetav\in\R^{\dimp}} \Lab(\thetav
,\thetav^{*}) $
can be called the \emph{upper excess} %\emph{value of the upper
%problem},
while
$ \sup_{\thetav\in\Thetas(\rups)} \Lam(\thetav,\thetav^{*}) $ is
the \emph{lower excess}.
The quadratic (in $ \thetav$) structure of the functions $ \Lab
(\thetav,\thetav^{*}) $
and $ \Lam(\thetav,\thetav^{*}) $ enables us to explicitly solve the
problem of
maximizing the corresponding function w.r.t. $ \thetav$.
%The quality of the quadratic approximation \eqref{LttbLtt} can be
%naturally measured
%by the difference between excesses of the upper and lower
%approximating problems.
%The approximating quality is sufficiently good if this difference is
%smaller in order
%than the value itself.
%

%le3.8 #&#
\begin{lemma}
\label{LLbreveloc}
It holds
%
%e3.23 #&#
\begin{equation}
\sup_{\thetav\in\R^{\dimp}} \Lab\bigl(\thetav,\thetav^{*}\bigr) = \| \xivb
\|^{2}/2 . \label{supLat}
\end{equation}
On the random set $ \{ \| \VPc\DPm^{-1} \xivm\| \le\rr\} $, it
also holds
%
%e3.24 #&#
\[
\sup_{\thetav\in\Thetas(\rr)} \Lam(\thetav,\thetav) = \| \xivm\|^{2}/2 .
% \le
% \| \xiv\|^{2}/2 .
\label{supLatm}
\]
\end{lemma}
\begin{pf}
The unconstrained maximum of the quadratic form $ \Lab(\thetav
,\thetav^{*}) $ w.r.t.~$ \thetav$ is attained at
$ \tilde{\thetav}_{\rd} = \DPb^{-1} \xivb= \DPb^{-2} \nabla
L(\thetav^{*}), $
yielding the expression \eqref{supLat}.
The lower excess is computed similarly.
\end{pf}

Our next step is in bounding the difference
$ \| \xivb\|^{2} - \| \xivm\|^{2} $.
It can be decomposed as
%
%e3.25 #&#
\[
%[c]
\| \xivb\|^{2} - \| \xivm\|^{2} = \| \xivb
\|^{2} - \| \xiv\|^{2} + \| \xiv\|^{2} - \| \xivm
\|^{2} \label{xivbvvm}
\]
with $ \xiv= \DPc^{-1} \nabla L(\thetav^{*}) $.
If the values $ \rddelta, \rdomega$ are small, then the difference
$ \| \xivb\|^{2} - \| \xivm\|^{2} $ is automatically smaller than
$ \| \xiv\|^{2} $.
%Under the identifiability condition \( (\AssId) \), this claim can be
%justified by
%the following simple lemma.

%le3.9 #&#
\begin{lemma}
\label{Lxivgap}
Suppose $ (\AssId) $ and let
$ \tau_{\rd} \eqdef\rddelta+ \rdomega\fis^{2} < 1 $.
Then
%
%e3.26 #&#
%e3.27 #&#
\begin{eqnarray}\label{alprdpa}
%[c]
\DPb^{2}& \ge& (1 - \tau_{\rd}) \DPc^{2} ,\qquad
\DPm^{2} \le (1 + \tau_{\rd}) \DPc^{2} ,
\nonumber
\\[-8pt]
\\[-8pt]
\nonumber
\bigl\| \Id_{\dimp} - \DPb\DPm^{-2} \DPb\bigr\|_{\infty}& \le&
\alp_{\rd} \eqdef \frac{2 \tau_{\rd}}{1 - \tau_{\rd}^{2}} .
\end{eqnarray}
Moreover,
%
%e3.28 #&#
%e3.29 #&#
\begin{eqnarray*}
%[c]
\| \xivb\|^{2} - \| \xiv\|^{2} &\le&
\frac{\tau_{\rd}}{1 - \tau_{\rd}} \| \xiv\|^{2}, \qquad \| \xiv\|^{2} - \| \xivm
\|^{2} \le \frac{\tau_{\rd}}{1 + \tau_{\rd}} \| \xiv\|^{2},
\\
\| \xivb\|^{2} - \| \xivm\|^{2}& \le &\alp_{\rd} \|
\xiv\|^{2} . \label{xivbccmpa}
\end{eqnarray*}
\end{lemma}

Our final step is in showing
that under $ (E D_{0}) $, the norm $ \| \xiv\| $ behaves essentially
as a norm of a Gaussian vector with the same covariance matrix.
Define for $ \BB\eqdef\DPc^{-1} \VPc^{2} \DPc^{-1} $
%
%e3.30 #&#
\[
%[c]
\dimA \eqdef \tr ( \BB ) ,\qquad \vA^{2} \eqdef 2 \tr\bigl(
\BB^{2}\bigr),\qquad \lambda_{0} \eqdef \| \BB\|_{\infty} =
\lambda_{\max} ( \BB ). \label{BBrdd}
\]
Under the identifiability condition $ (\AssId) $, one can bound
%
%e3.31 #&#
\[
%[c]
\BB^{2} \le \fis^{2} \Id_{\dimp},\qquad \dimA \le
\fis^{2} \dimp, \qquad\vA^{2} \le 2 \fis^{4} \dimp,\qquad
\lambda_{0} \le \fis^{2}. \label{BBrdreg}
\]
Similarly to the previous result, we assume that the constant
$ \gmb$ from condition $ (E D_{0}) $ is sufficiently
large, namely, $ \gm^{2} \ge2 \dimA$.
%The other case only changes the constants in the inequalities.
Define $ \muc= 2/3 $
%Moreover, with the constant \( \gmb\) from \( (E D_{0}) \),
and %\( \muc= \gm^{2}/(\dimA+ \gm^{2}) \wedge2/3 \), and
%
%e3.32 #&#
%e3.33 #&#
%e3.34 #&#
\begin{eqnarray*}
\label{yycgmcxxcAd} \yyc^{2} & \eqdef& \gmb^{2}/
\muc^{2} - \dimA/\muc,
\\
\gmc & \eqdef& \muc\yyc = \sqrt{\gm^{2} - \muc\dimA} ,
\\
2\xxc & \eqdef& \muc\yyc^{2} + \log\det \bigl( \Id_{\dimp} -
\muc\BB^{2} / \lambda_{0} \bigr) .
\end{eqnarray*}
It is easy to see that $ \yyc^{2} \ge3 \gmb^{2}/2 $ and
$ \gmc\ge\sqrt{2/3}   \gmb$.

%th3.10 #&#
\begin{theorem}
\label{LLbrevelocm}
Let $ (E D_{0}) $ hold with $ \nunu= 1 $ and
$ \gmb^{2} \ge2 \dimA$.
Then
$ \E\| \xiv\|^{2} \le\dimA$, and for each $ \xx\le\xxc$
%
%e3.35 #&#
\begin{equation}
\P \bigl( \| \xiv\|^{2}/\lambda_{0} \ge\zz(\xx,\BB) \bigr)
\le 2 \ex^{-\xx} + 8.4 \ex^{-\xxc} , \label{PxivbzzBB}
\end{equation}
where $ \zz(\xx,\BB) $ is defined by
%
%e3.36 #&#
\[
\label{PzzxxpB} \zz(\xx,\BB) \eqdef \cases{ \dimA+ 2 \vA\xx^{1/2}, &\quad  $
\xx\le\vA/18 ,$ \vspace*{2pt}
\cr
\dimA+ 6 \xx, &\quad  $\vA/18 < \xx\le\xxc.$ }
% \bigl| \yyc+ 2 (\xx- \xxc)/\gmc\bigr|^{2}, & \xx> \xxc.
\label{zzxxppdBl}
\]
Moreover, for $ \xx> \xxc$, it holds with
$ \zz(\xx,\BB) =  | \yyc+ 2 (\xx- \xxc)/\gmc |^{2} $
%
%e3.37 #&#
\[
\P \bigl( \| \xiv\|^{2}/\lambda_{0} \ge\zz(\xx,\BB) \bigr)
\le 8.4 \ex^{-\xx} . \label{PxivbzzBBc}
\]
\end{theorem}

\begin{pf}
It follows from condition $ (E D_{0}) $ that
%
%e3.38 #&#
%e3.39 #&#
\begin{eqnarray*}
\label{Exiv2tr} \E\| \xiv\|^{2} &=& \E\tr\xiv\xiv^{\T}
\\
&=& \tr\DPc^{-1} \bigl[ \E\nabla L\bigl(\thetav^{*}\bigr)
\bigl\{ \nabla L\bigl(\thetav^{*}\bigr) \bigr\}^{\T} \bigr]
\DPc^{-1} = \tr \bigl[ \DPc^{-2} \Var \bigl\{ \nabla L\bigl(
\thetav^{*}\bigr) \bigr\} \bigr]
\end{eqnarray*}
and $ (E D_{0}) $ implies
$ \gammav^{\T} \Var \{ \nabla L(\thetav^{*})  \} \gammav
\le\gammav^{\T} \VPc^{2} \gammav$ and, thus,
$ \E\| \xiv\|^{2} \le\dimA$.
The deviation bound \eqref{PxivbzzBB} is proved in Corollary~2.5 %
of the supplement [\citet{supp}].
\end{pf}

%re3 #&#
\begin{remark}
\label{Rxxc}
This small remark concerns the term $ 8.4 \ex^{-\xxc} $ in the
probability bound
\eqref{PxivbzzBB}.
As already mentioned, this bound implicitly assumes that the constant
$ \gmb$ is
large (usually $ \gmb\asymp\nsize^{1/2} $).
Then $ \xxc\asymp\gmb^{2} \asymp\nsize$ is large as well.
So, $ \ex^{-\xxc} $ is very small and asymptotically negligible.
Below we often ignore this term.
For $ \xx\le\xxc$, we can use $ \zz(\xx,\BB) = \dimA+ 6 \xx$.
\end{remark}

%re4 #&#
\begin{remark}
\label{RLLbrevelocm}
The exponential bound of Theorem~\ref{LLbrevelocm} helps to describe
the critical value of $ \zz$ ensuring
a prescribed deviation probability
$ \P ( \| \xiv\|^{2} \ge\zz ) $.
Namely, this probability starts to gradually decrease
when $ \zz$ grows over $ \lambda_{0} \dimA$.
In particular, this helps to answer a very important question about a
proper choice
of the critical value $ \zz$ providing the prescribed covering probability,
or of the value $ z $ ensuring the dominating concentration probability
$ \P ( \| \DPb(\tilde{\thetav} - \thetav^{*}) \| \le z  )
$.

The definition of the set $ \LCS_{\rd}(\rr) $ from \eqref
{LCSrdrr} involves
the event $ \{ \| \VPc\DPm^{-1} \xivm\| > \rr\} $.
Under $ (\AssId) $, it is included in the set
$ \{ \| \xivm\| > (1 + \alp_{\rd})^{-1} \fis^{-1} \rr\} $ [see
\eqref{alprdpa}], and
its probability is of order $ \ex^{-\xx} $ for
$ \rr^{2} \ge C (\xx+ \dimp) $ with a fixed $ C > 0 $.
\end{remark}

By Theorem~\ref{Tdefbound}, one can use
$ \max \{ \errb(\rr), \errm(\rr)  \} \le\rdomega
\zzQ(\xx,\entrlb) $
on a set of probability at least $ 1 - \ex^{-\xx} $.
Further, $ \| \xiv\|^{2} / \lambda_{0} \le\zz(\xx,\BB) $ with a
probability of
order $ 1 - \ex^{-\xx} $;
see \eqref{PxivbzzBB}.
Putting together the obtained bounds yields for the spread
$ \spread_{\rd}(\rr) $ with a probability
about $ 1 - 4 \ex^{-\xx} $
%
%e3.40 #&#
\[
\spread_{\rd}(\rr) \le 2 \rdomega \zzQ(\xx,\entrlb) +
\alp_{\rd} \lambda_{0} \zz(\xx,\BB) . \label{defiregc}
\]
The results obtained in Section~\ref{Slocinfr} are sharp and
meaningful if the
spread $ \spread_{\rd}(\rr) $ is smaller in order than
the value $ \E\| \xiv\|^{2} $.
Theorem~\ref{LLbrevelocm} states that $ \| \xiv\|^{2} $
does not significantly deviate over its expected value
$ \dimA\eqdef\E\| \xiv\|^{2} $ which is our leading term.
We know that
$ \zzQ(\xx,\entrlb) \approx\entrlb+ \xx= \cdimb\dimp+ \xx$ if
$ \xx$ is not too large.
Also, $ \zz(\xx,\BB) \le\dimA+ 6 \xx$, where $ \dimA$ is of
order $ \dimp$
due to $ (\AssId) $.
Summarizing the above discussion yields that the local results apply if
the regularity
condition $ (\AssId) $ holds and the values
$ \rdomega$ and $ \alp_{\rd} $ or, equivalently,
$ \rhor(\rr), \rddelta(\rr) $ are small.
In Section~\ref{Chgexamples} we show for the i.i.d. example that
$ \rhor(\rr) \asymp\sqrt{\rr^{2}/\nsize} $ and similarly for $
\rddelta(\rr) $.

%The bound \eqref{taurhodel} can be further specified in the so called
%under the condition
% \VPc\le\fis\DPc  .
%If the parametric assumption is correct, that is, \( \P= \P_{
%a regular parametric family,
%then the both matrices coincide with the total Fisher information
%matrix, and the
%regularity condition is fulfilled automatically with \( \fis= 1 \).
%Otherwise, the regularity means that the local variability of the
%process
%significantly larger than the
%local information measured by the matrix \( \DPc\).
%
%Suppose \eqref{VPcfisDPc} and let
%Then
% \DPb^{2}
% \ge
% (1 - \tau_{\rd}) \DPc^{2} ,
%
% \alp_{\rd}
% =
% \| \Id_{\dimp} - \DPb\DPm^{-2} \DPb\|_{\infty}
% \le
% \frac{2 \tau_{\rd}}{1 - \tau_{\rd}}.
%Moreover, \( \BB^{2} = \DPb^{-1} \VPc^{2} \DPb^{-1} \) satisfies with
% \BB^{2}
% \le
% \norma_{\rd} \Id_{\dimp},
%
% \dimA
% \le
% \norma_{\rd} \dimp,
%
% \vA^{2}
% \le
% 2 \norma_{\rd}^{2} \dimp,
%
% \lambda_{0}
% \le
% \norma_{\rd}.
%
%The results follow directly from the definition of \( \DPb\) and
%
%In particular,
%the matrices \( \DPb\) and \( \DPm\) are close to each other if \(
%and \( \rdomega\fis^{2} \) are small.
%%
%So, all we need in the regular case, is a large deviation bound for
%the probability \( \P\bigl( \tilde{\thetav} \notin\Thetas(\rr)
%and that the quantities \( \rhor(\rr) \) and \( \rddelta(\rr) \) are
%small.
%

%s3.4 #&#
\subsection{\texorpdfstring{Proof of Corollaries~\protect\ref{TquadWilks} and \protect\ref{Tquadexploc}}
{Proof of Corollaries 3.3 and 3.4}}
The bound \eqref{LttbLtt} together with Lem\-ma~\ref{LLbreveloc} yield on
$ \LCS_{\rd}(\rr) $
%
%e3.41 #&#
%e3.42 #&#
\begin{eqnarray}
\label{LttRRLLa} L\bigl(\tilde{\thetav},\thetav^{*}\bigr) &=&
\sup_{\thetav\in\Thetas(\rr)} L\bigl(\thetav,\thetav^{*}\bigr)
\nonumber
\\[-8pt]
\\[-8pt]
\nonumber
& \ge& \sup_{\thetav\in\Thetas(\rr)} \Lam\bigl(\thetav,\thetav^{*}\bigr) - \errm
(\rr) = \| \xivm\|^{2}/2 - \errm(\rr) .
\end{eqnarray}
Similarly,
%
%e3.43 #&#
\[
L\bigl(\tilde{\thetav},\thetav^{*}\bigr) \le \sup_{\thetav\in\Thetas(\rr)} \Lab
\bigl(\thetav,\thetav^{*}\bigr) + \errb (\rr) \le \| \xivb
\|^{2}/2 + \errb(\rr), \label{LttRRLLb}
\]
yielding \eqref{VPctttnL}.
For getting \eqref{DPbttt}, we again apply the inequality
$ L(\thetav,\thetav^{*}) \le\Lab(\thetav,\thetav^{*}) + \errb
(\rr)
$ from Theorem~\ref{TapproxLL}
for $ \thetav$ equal to $ \tilde{\thetav} $.
With $ \xivb= \DPb^{-1} \nabla L(\thetav^{*}) $ and
$ \uv_{\rd} \eqdef\DPb(\tilde{\thetav} - \thetav^{*}) $,
this gives
%
%e3.44 #&#
\[
L\bigl(\tilde{\thetav},\thetav^{*}\bigr) - \xivb^{\T}
\uv_{\rd} + \| \uv_{\rd} \|^{2}/2 \le \errb(\rr) .
\label{ttD2rK}
\]
Therefore, by \eqref{LttRRLLa},
%
%e3.45 #&#
\[
%[c]
\| \xivm\|^{2}/2 - \errm(\rr) - \xivb^{\T}
\uv_{\rd} + \| \uv_{\rd} \|^{2}/2 \le \errb(\rr)
\label{f12xivtxiv}
\]
or, equivalently,
%
%e3.46 #&#
\[
%[c]
\| \xivb\|^{2}/2 - \xivb^{\T} \uv_{\rd} +
\| \uv_{\rd} \|^{2}/2 \le \errb(\rr) + \errm(\rr) + \bigl( \|
\xivb\|^{2} - \| \xivm\|^{2} \bigr)/2 \label{f12xivtxiv0}
\]
and the definition of $ \spread_{\rd}(\rr) $ implies
$  \| \uv_{\rd} - \xivb \|^{2} \le2 \spread_{\rd}(\rr
) $.

%================================ LD ================================
%s4 #&#
\section{Upper function approach and concentration of the qMLE}
\label{Chgparam}

%The whole proposed approach relies implicitly on the two groups of
%assumptions:
%global and local.
%These assumptions are linked to each other by the value \( \rr\).
%From one side, the \emph{global} assumptions listed
%in Section~\ref{Sglobalcondition} and used in Theorem~\ref{Thittinggl}
%and
%its corollaries below in Section~\ref{Chgparam} should ensure a
%sensible bound for
%the deviation probability
%This particularly requires that \( \rr\) is sufficiently large.
%In the contrary, the \emph{local} conditions are based on the
%assumption
%that the local set \( \Thetas(\rr) \) is sufficiently small to
%guarantee that
%errors of approximation \( \errb(\rr) \) and \( \errm(\rr) \) and
%the spread \( \spread_{\rd}(\rr) \) from \eqref{taurhodel} are small
%as well
%in a probabilistic sense.

A very important step in the analysis of the qMLE $ \tilde{\thetav}
$
is \emph{localization}.
This property means that $ \tilde{\thetav} $ concentrates in a
small vicinity
of the central point $ \thetav^{*}$.
This section states such a concentration bound under the global
conditions of
Section~\ref{Scondgllo}.
Given $ \rups$,
the deviation bound describes the probability
$ \P ( \tilde{\thetav} \notin\Thetas(\rups)  ) $ that
$ \tilde{\thetav} $ does not belong to the local vicinity $
\Thetas(\rups) $ of
$ \Theta$.
The question of interest is to check a possibility of selecting $
\rups$ in a way
that the local bracketing result and the deviation bound apply simultaneously;
see the discussion at the end of the section.

Below we suppose that a sufficiently large constant $ \xx$ is fixed
to specify
the accepted level be of order $ \ex^{-\xx} $ for this deviation
probability.
All the constructions below depend upon this constant.
We do not indicate it explicitly for ease of notation.

The key step in this large deviation bound is made in terms of an
\emph{upper function} for the process
$ L(\thetav,\thetav^{*}) \eqdef L(\thetav) - L(\thetav^{*}) $.
Namely, $ \pnnd(\thetav) $ is a deterministic \emph{upper
function} if it holds
with a high probability:
%
%e4.1 #&#
\begin{equation}
%[c]
\sup_{\thetav\in\Theta} \bigl\{ L\bigl(\thetav,\thetav^{*}
\bigr) + \pnnd(\thetav) \bigr\} \le0. \label{PLttZzz0}
\end{equation}
Such bounds are usually called for in the analysis of the posterior
measure in
the Bayes approach.
Below we % suppose that such a function \( \pnnd(\cdot) \) is fixed and
present sufficient conditions ensuring \eqref{PLttZzz0}.
Now we explain how \eqref{PLttZzz0}
can be used for describing \emph{concentration sets} for~$\tilde{\thetav}$.

%le4.1 #&#
\begin{lemma}
\label{CThittinggl}
Let $ \pnnd(\thetav) $ be an upper function in the sense
%
%e4.2 #&#
\begin{equation}
%[c]
\P \Bigl( \sup_{\thetav\in\Theta} \bigl\{ L\bigl(\thetav,
\thetav^{*}\bigr) + \pnnd(\thetav) \bigr\} \ge 0 \Bigr) \le
\ex^{-\xx} \label{hitprobxxgl}
\end{equation}
for $ \xx> 0 $.
Given a subset $ \Thetas\subset\Theta$ with $ \thetav^{*}\in
\Thetas$,
the condition $ \pnnd(\thetav) \ge0 $ for
$ \thetav\notin\Thetas$ ensures
%
%e4.3 #&#
\[
%[c]
\P ( \tilde{\thetav} \notin\Thetas ) \le \ex^{-\xx} .
\label{PCAxxgl}
\]
\end{lemma}

\begin{pf}
If $ \Theta^{\circ} $ is a subset of $ \Theta$ not containing $
\thetav^{*}$,
then the event $ \tilde{\thetav} \in\Theta^{\circ} $ is only
possible if
$ \sup_{\thetav\in\Theta^{\circ}} L(\thetav,\thetav^{*}) \ge0 $,
because $ L(\thetav^{*},\thetav^{*}) \equiv0 $.
%This yields the result.
\end{pf}

A possible way of checking the condition \eqref{hitprobxxgl} is based
on a lower
quadratic bound for the negative expectation
$ - \E L(\thetav,\thetav^{*}) \ge\gmi(\rr) \| \VPc(\thetav-
\thetav^{*}) \|^{2}/2 $
in the sense of condition $ (\cc{L}\rr) $ from Section~\ref
{Sglobalcondition}.
We present two different results.
The first one assumes that the values $ \gmi(\rr) $ can be fixed
universally for
all $ \rr\ge\rups$.
%
%th4.2 #&#
\begin{theorem}
\label{CThittingglrc}
Suppose $ (E\rr) $ and $ (\cc{L}\rr) $ with $ \gmi(\rr)
\equiv\gmi$.
Let, for $ \rr\ge\rups$,
%
%e4.4 #&#
%e4.5 #&#
\begin{eqnarray}
\label{cgmi1rrc} 1 + \sqrt{\xx+ \entrlb} & \le& 3 \nunu^{2} \gm(\rr)/
\gmi,
\\
6 \nunu\sqrt{\xx+ \entrlb} & \le& \rr\gmi, \label{cgmi2rrc}
\end{eqnarray}
with $ \xx+ \entrlb\ge2.5 $ and $ \entrlb= \cdimb\dimp$. Then
%
%e4.6 #&#
\begin{equation}
%[c]
\P \bigl( \tilde{\thetav} \notin\Thetas(\rups) \bigr) \le
\ex^{-\xx} . \label{PnotinTsruc}
\end{equation}
\end{theorem}

\begin{pf}
The result follows from Theorem~2.8 of the supplement [\citet{supp}] %\ref{Thitting}
with $ \mubc= \frac{\gmi}{3 \nunu} $,
$ \pen(\mubc) \equiv0 $,
$ \UP(\thetav) = L(\thetav) - \E L(\thetav) $ and
$ \Ldrift(\thetav,\thetav^{*}) = \break- \E L(\thetav,\thetav^{*})
\ge\frac{\gmi}{2} \| \VPc(\thetav- \thetav^{*}) \|^{2} $.
\end{pf}

%re5 #&#
\begin{remark}
\label{RCThittingglrc}
The bound \eqref{PnotinTsruc} requires only two conditions.
Condition \eqref{cgmi1rrc} means that the value $ \gm(\rr) $ from condition
$ (E\rr) $ fulfills
$ \gm^{2}(\rr) \ge C (\xx+ \dimp) $, that is, we need a qualified
rate in the
exponential moment conditions.
This is similar to requiring finite polynomial moments for the score function.
Condition \eqref{cgmi2rrc} requires that $ \rr$ exceeds some fixed
value, namely,
\mbox{$ \rr^{2} \ge C (\xx+ \dimp) $}.
This bound is helpful for fixing the value $ \rups$ providing a
sensible deviation
probability bound.
\end{remark}

If $ \gmi(\rr) $ decreases with $ \rr$, the result is a bit more
involved.
The key requirement is that $ \gmi(\rr) $ decreases not too fast,
so that the product $ \rr\gmi(\rr) $ grows to infinity with~$ \rr
$.
The idea is to include the complement of the central set $ \Thetas$
in $ \Theta$
in the union of the growing
sets $ \Thetas(\rr_{k}) $ with $ \gmi(\rr_{k}) \ge\gmi(\rups)
2^{-k} $,
and then apply Theorem~\ref{CThittingglrc} for each $ \Thetas(\rr_{k}) $.

%th4.3 #&#
\begin{theorem}
\label{CThittingglr}
Suppose $ (E\rr) $ and $ (\cc{L}\rr) $.
Let $ \rr_{k} $ be such that $ \gmi(\rr_{k}) \ge\gmi(\rups)
2^{-k} $ for
$ k \ge1 $.
If the conditions
%
%e4.7 #&#
%e4.8 #&#
\begin{eqnarray*}
\label{cgmi2rr} 1 + \sqrt{\xx+ \entrlb+ c k} & \le& 3 \nunu^{2} \gm(
\rr_{k})/\gmi(\rr_{k}) ,
\\
6 \nunu\sqrt{\xx+ \entrlb+ c k} & \le& \rr_{k} \gmi(
\rr_{k}) ,
\end{eqnarray*}
are fulfilled for $ c = \log(2) $, then it holds
%
%e4.9 #&#
\[
%[c]
\P \bigl( \tilde{\thetav} \notin\Thetas(\rups) \bigr) \le
\ex^{-\xx} . \label{ttnotinTsrs}
\]
\end{theorem}

\begin{pf}
The result \eqref{PnotinTsruc} is applied to each set $ \Thetas(\rr_{k}) $ and
$ \xx_{k} = \xx+ c k $.
This yields
%
%e4.10 #&#
\[
%[c]
\P \bigl( \tilde{\thetav} \notin\Thetas(\rups) \bigr) \le \sum
_{k \ge1} \P \bigl( \tilde{\thetav} \notin\Thetas(
\rr_{k}) \bigr) \le \sum_{k \ge1}
\ex^{-\xx- c k} = \ex^{-\xx} \label{PnotinTsrucm}
\]
as required.
\end{pf}

%This result will be further specified in Section~\ref{SqMLEiid} for
%the i.i.d. model.

%re6 #&#
\begin{remark}
\label{RLLDLD}
Here we briefly discuss the very important question: how one can fix
the value
$ \rups$ ensuring the bracketing result in the local set $ \Thetas
(\rups) $ and
a small probability of the related set $ \LCS_{\rd}(\rr) $ from
\eqref{LCSrdrr}?
The event $ \{ \| \VPc\DPm^{-1} \xivm\| > \rr\} $ requires
$ \rr^{2} \ge C (\xx+ \dimp) $.
Further, we inspect
the deviation bound for the complement $ \Theta\setminus\Thetas
(\rups) $.
For simplicity, assume $ (\cc{L}\rr) $ with $ \gmi(\rr) \equiv
\gmi$.
Then the condition~\eqref{cgmi2rrc} of Theorem~\ref{CThittingglrc}
requires that
%
%e4.11 #&#
\begin{equation}
%[c]
\rups^{2} \ge C \gmi^{-2} (\xx+ \dimp) .
\label{rups2gmixxd}
\end{equation}
In words, the squared radius $ \rups^{2} $ should be at least of
order $ \dimp$.
The other condition~\eqref{cgmi1rrc} of Theorem~\ref{CThittingglrc}
is technical and
only requires that $ \gm(\rr) $ is sufficiently large, while
the local results only require that $ \rddelta(\rr) $ and $
\rdomega(\rr) $ are
small for such $ \rr$.
In the asymptotic setup one can typically bring these conditions together.
Section~\ref{Chgexamples} provides further discussion for the i.i.d. setup.
\end{remark}

%================================ exc ================================

%s5 #&#
\section{Examples}
\label{Chgexamples}
\label{Sexpex}

The model with independent identically distributed (i.i.d.) observations
is one of the most popular setups in statistical literature and
in statistical applications.
The essential and the most developed part of the statistical theory is designed
for the i.i.d. modeling.
Especially, the classical asymptotic parametric theory is almost
complete including
asymptotic root-n normality and efficiency of the MLE
and Bayes estimators under rather mild assumptions; see, for example,
Chapters 2 and 3 in
\citet{IH1981}.
So, the i.i.d. model can naturally serve as a benchmark for any
extension of the statistical
theory: being applied to the i.i.d. setup, the new approach should lead to
essentially the same conclusions as in the classical theory.
Similar reasons apply to the regression model and its extensions.
Below we try to demonstrate that the proposed nonasymptotic viewpoint
is able to reproduce the existing brilliant and well-established
results of the
classical parametric theory. % for the i.i.d. and regression setup.
Surprisingly,
the majority of classical efficiency results can be easily derived from
the obtained general nonasymptotic bounds.

The next question is whether there is any added value or benefits of
the new approach
being restricted to the i.i.d. situation relative to the classical one.
Two important
issues have been already mentioned: the new approach applies to the
situation with
finite samples and survives under model misspecification. One more
important question is
whether the obtained results remain applicable and informative if the
dimension of the
parameter space is high---this is one of the main challenges in modern
statistics.
We show that the dimensionality $ \dimp$ naturally appears in the
risk bounds and the
results apply as long as the sample size exceeds in order of this value
$ \dimp$.
%The methods of reducing dimensionality or complexity of the model are
%discussed in
%Chapters~\ref{Srough} and~\ref{Sposterior}.
All these questions are addressed in Section~\ref{SqMLEiid} for the
i.i.d. setup;
Section~\ref{SexpexGLM} focuses on
generalized linear modeling, while Section~\ref{Slinmedian} discusses
linear median
regression.

%s5.1 #&#
\subsection{Quasi-MLE in an i.i.d. model}
\label{SqMLEiid}
An i.i.d. parametric model means that the observations
$ \Yv= (Y_{1},\ldots,Y_{\nsize}) $ are independent identically
distributed from a
distribution $ P $ which belongs to a given parametric family
$ (P_{\thetav}, \thetav\in\Theta) $ on the observation space $
\YY_{1} $.
Each $ \thetav\in\Theta$
clearly yields the product data distribution
$ \P_{\thetav} = P_{\thetav}^{\otimes\nsize} $
on the product space $ \YY= \YY_{1}^{\nsize} $.
This section illustrates how the obtained general results can be
applied to this
type of modeling under possible model misspecification.
Different types of misspecification can be considered. Each of the
assumptions, namely,
data independence, identical distribution and parametric form of the marginal
distribution can be violated.
To be specific, we assume the observations $ Y_{i} $ independent and
identically distributed.
However, we admit that
%two types of misspecification: the \( Y_{i} \)'s could be
%differently distributed and
the distribution of each $ Y_{i} $ does not necessarily
belong to the parametric family $ (P_{\thetav}) $.
The case of nonidentically distributed observations can be done
similarly at the cost
of more complicated notation.

In what follows the parametric family $ (P_{\thetav}) $ is supposed
to be dominated
by a measure $ \mu_{0}$, and each density $ p(y,\thetav) =
dP_{\thetav}/d\mu_{0}(y) $
is two times continuously differentiable in $ \thetav$ for all $ y
$.
Denote $ \ell(y,\thetav) = \log p(y,\thetav) $.
The parametric assumption $ Y_{i} \sim P_{\thetav^{*}} \in(P_{\thetav
}) $ leads to
the log-likelihood
%
%e5.1 #&#
\[
%[c]
L(\thetav) = \sum\ell(Y_{i},\thetav) ,
\label{Ltiid}
\]
where the summation is taken over $ i=1,\ldots,\nsize$.
The quasi-MLE $ \tilde{\thetav} $ maximizes this sum over $
\thetav\in\Theta
$:
%
%e5.2 #&#
\[
%[c]
\tilde{\thetav} \eqdef \mathop{\argmax}_{\thetav\in\Theta} L(\thetav) =
\mathop{\argmax}_{\thetav\in\Theta} \sum\ell(Y_{i},\thetav). \label{tttiid}
\]
The target of estimation $ \thetav^{*}$ maximizes the expectation of $
L(\thetav)
$:
%
%e5.3 #&#
\[
%[c]
\thetav^{*} \eqdef \mathop{\argmax}_{\thetav\in\Theta} \E L(\thetav) =
\mathop{\argmax}_{\thetav\in\Theta} \E\ell(Y_{1},\thetav). \label{tsiid}
\]
Let $ \zeta_{i}(\thetav) \eqdef\ell(Y_{i},\thetav) - \E\ell
(Y_{i},\thetav) $.
Then $ \zeta(\thetav) = \sum\zeta_{i}(\thetav) $.
The equation\break
$ \nabla\E L(\thetav^{*}) = 0 $ implies
%
%e5.4 #&#
\begin{equation}
%[c]
\nabla\zeta\bigl(\thetav^{*}\bigr) = \sum\nabla
\zeta_{i}\bigl(\thetav^{*}\bigr) = \sum\nabla
\ell_{i}\bigl(\thetav^{*}\bigr) . \label{nztiid}
\end{equation}

I.i.d. structure of the $ Y_{i} $'s allows to rewrite the local conditions
$ (E\rr) $, $ (E D_{0}) $, $ (E D_{1}) $, and $ (\LL_{0}) $
, and
$ (\AssId) $
in terms of the marginal distribution.
%In the following conditions the index \( i \) runs from \( 1 \) to \(
%
\begin{longlist}[$({ed_{0}})$]
% \emph{ There exist some constants \( \nunu\),
% %continuous symmetric matrix function \( V(\thetav) \) for \( \thetav
% and \( \gmiid> 0 \), and %for each \( i \),
% a positive symmetric \( \dimp\times\dimp\) matrix \( \vp\),
% such that for all \( |\lambda| \le\gmiid\)
% }
% \sup_{\gammav\in\R^{\dimp}} \sup_{\thetav\in\Theta}
% \log\E\exp\biggl\{
% \lambda\frac{\gammav^{\T} \nabla\zeta_{i}(\thetav)}
% {\| \vp\gammav\|}
% \biggr\} \le
% \nunu^{2} \lambda^{2} / 2.

\item[$({ed_{0}})$]
\textit{There exists a positively definite symmetric matrix $ \mathbf
{v}_{0}$,
such that for all $ |\lambda| \le\gmiid$}
%
%e5.5 #&#
\[
%[c]
\label{expzetaciid} \sup_{\gammav\in\R^{\dimp}} \log\E\exp \biggl\{ \lambda
\frac{\gammav^{\T} \nabla\zeta_{1}(\thetav^{*})} {
\| \mathbf{v}_{0}\gammav\|} \biggr\} \le \nunu^{2} \lambda^{2} / 2.
\]
\end{longlist}
A natural candidate on $ \mathbf{v}_{0}^{2} $ is given by the
variance of the gradient
$ \nabla\ell(Y_{1},\thetav^{*}) $, that is,
$ \mathbf{v}_{0}^{2} = \Var \{ \nabla\ell(Y_{1},\thetav^{*})
\}
= \Var \{ \nabla\zeta_{1}(\thetav^{*})  \} $.

Next consider the local sets
%
%e5.6 #&#
\[
%[c]
\Thetasi(\rri) = % \{ \thetav: \| \VPc(\thetav- \thetavs) \| \le\rr\}
% =
\bigl\{ \thetav\dvtx \bigl\|
\mathbf{v}_{0}\bigl(\thetav- \thetav^{*}\bigr) \bigr\| \le\rri
\bigr\} . \label{Theta0riid}
\]
In view of $ \VPc^{2} = \nsize\mathbf{v}_{0}^{2} $, it holds
$ \Thetas(\rr) = \Thetasi(\rri) $ with $ \rr^{2} = \nsize\rri^{2} $.

Below we distinguish between local conditions for $ \rri\le\rris$
and the global conditions for all $ \rri> 0 $, where
$ \rris$ is some fixed value.

The local smoothness conditions $ (E D_{1}) $ and $ (\LL_{0}) $
require to specify the functions $ \rddelta(\rr) $ and $ \rdomega
(\rr) $
for $ \rr\le\rups$ where $ \rups^{2} = \nsize\rris^{2} $.
%We specify these conditions assuming that the considered parametric
%family is
%sufficiently regular.
If the log-likelihood function $ \ell(y,\thetav) $ is sufficiently
smooth in
$ \thetav$, these functions can be selected proportional to
$ \rri= \rr/\nsize^{1/2} $.

\begin{longlist}[$({ed_{1}})$]
\item[$({ed_{1}})$]
\textit{There are constants $ \rhor^{*}> 0 $ and $ \gmiid> 0 $ such
that for
each $ \rri\le\rris$ and $ |\lambda| \le\gmiid$}
%
%e5.7 #&#
\[
%[c]
\label{expzetac0iid} \sup_{\gammav\in\R^{\dimp}} \sup_{\thetav\in\Thetasi(\rri)} \log\E
\exp \biggl\{ \lambda \frac{\gammav^{\T}  [
\nabla\zeta_{1}(\thetav) - \nabla\zeta_{1}(\thetav^{*})
]} {
\rhor^{*}  \rri  \| \mathbf{v}_{0}\gammav\|} \biggr\} \le \nunu^{2}
\lambda^{2} / 2 .
\]
\end{longlist}

Further, we restate the local identifiability condition $ (\LL_{0})
$ in terms of
the expected value
$ \kullbi(\thetav,\thetav^{*})
\eqdef- \E \{ \ell(Y_{i},\thetav) - \ell(Y_{i},\thetav^{*})
\} $
for each $ i $.
We suppose that $ \kullbi(\thetav,\thetav^{*}) $ is two times
differentiable w.r.t. $ \thetav$.
The definition of $ \thetav^{*}$ implies $ \nabla\E\ell
(Y_{i},\thetav^{*}) = 0 $.
Define also the matrix
$ \IFc= - \nabla^{2} \E\ell(Y_{i},\thetav^{*}) $.
In the parametric case $ P = P_{\thetav^{*}} $,
$ \kullbi(\thetav,\thetav^{*}) $ is the \emph{Kullback--Leibler
divergence} between
$ P_{\thetav^{*}} $ and $ P_{\thetav} $, while the matrices $
\mathbf{v}_{0}^{2} = \IFc$
are equal to each other and coincide with the
\emph{Fisher information matrix} of the family $ (P_{\thetav}) $ at
$ \thetav^{*}$.

\begin{longlist}[$({\ell_{0}})$]
\item[$({\ell_{0}})$]
\textit{There is a constant $ \rddelta^{*}$ such that it holds for each
$ \rri\le\rris$
% it holds on \( \Thetasi(\rri) \)
}
%
%e5.8 #&#
\[
%[c]
\label{LmgfquadELiid} \sup_{\thetav\in\Thetasi(\rri)} \biggl| \frac{2 \kullbi(\thetav,\thetav^{*})
% - (\thetav- \thetavs)^{\T} \nabla\kullbi(\thetavs)
} {
(\thetav- \thetav^{*})^{\T} \IFc  (\thetav- \thetav^{*}) } - 1 \biggr|
\le \rddelta^{*}\rri.
\]
%

%The \emph{identifiability condition} in terms of
%the family \( (P_{\thetav}) \) reads as
%
\item[$(\iota)$]
There is a constant $ \fis> 0 $ such that $ \fis^{2} \IFc^{2} \ge
\mathbf{v}_{0}^{2} $.

%Now we restate the global conditions \( (E\rr) \) and \( (\cc{L}\rr)

%
\item[$(e\rri)$]
\textit{For each $ \rri> 0 $, there exists $ \gmiid(\rri) > 0 $,
such that for all $ |\lambda| \le\gmiid(\rri) $}
%
%e5.9 #&#
\[
%[c]
\label{expzetaciid} \sup_{\gammav\in\R^{\dimp}} \sup_{\thetav\in\Thetasi(\rri)} \log\E
\exp \biggl\{ \lambda\frac{\gammav^{\T} \nabla\zeta_{1}(\thetav)} {
\| \mathbf{v}_{0}\gammav\|} \biggr\} \le \nunu^{2}
\lambda^{2} / 2.
\]

\item[$(\ell\rri)$]
\textit{For each $ \rri> 0 $, there exists $ \gmi(\rri) > 0 $
such that}
%
%e5.10 #&#
\[
\sup_{\thetav\in\Theta\dvtx    \| \mathbf{v}_{0}(\thetav- \thetav^{*})
\| = \rri} \frac{\kullbi(\thetav,\thetav^{*})}{\| \mathbf{v}_{0}(\thetav-
\thetav^{*}) \|^{2}} \ge \gmi(\rri) , \label{xxentrttiid}\vadjust{\goodbreak}
\]
%
%for \( b \ge12 \nunu\) and \( b_{1} \ge\).
\end{longlist}

%le5.1 #&#
\begin{lemma}
\label{Lcondiid}
Let $ Y_{1},\ldots,Y_{\nsize} $ be i.i.d.
Then $ (e\rri) $, $ (ed_{0}) $, $ (ed_{1}) $, $ (\iota) $
and $ (\ell_{0}) $ imply
$ (E\rr) $, $ (E D_{0}) $, $ (E D_{1}) $, $ (\AssId) $ and
$ (\LL_{0}) $ with
$ \VPc^{2} = \nsize\mathbf{v}_{0}^{2} $,
$ \DPc^{2} = \nsize\IFc$,
$ \rhor(\rr) = \rhor^{*}\rr/\nsize^{1/2} $,
$ \rddelta(\rr) = \rddelta^{*}\rr/\nsize^{1/2} $,
%the same constant \( \nunu\),
and $ \gmb= \gmiid\sqrt{\nsize} $.
\end{lemma}
\begin{pf}
The identities
$ \VPc^{2} = \nsize\mathbf{v}_{0}^{2} $,
$ \DPc^{2} = \nsize\IFc$
follow from the i.i.d. structure of the observations $ Y_{i} $.
We briefly comment on condition $ (E\rr) $.
The use of the i.i.d. structure once again yields by \eqref{nztiid}
in view of $ \VPc^{2} = \nsize\mathbf{v}_{0}^{2} $
%
%e5.11 #&#
\[
%[c]
\log\E\exp \biggl\{ \lambda\frac{\gammav^{\T} \nabla\zeta(\thetav)}{\| \VPc\gammav
\|} \biggr\} = \nsize\E\exp
\biggl\{ \frac{\lambda}{\nsize^{1/2}} \frac{\gammav^{\T} \nabla\zeta_{1}(\thetav)}{\| \mathbf
{v}_{0}\gammav\|} \biggr\} \le \nunu^{2}
\lambda^{2}/2 \label{Elognsizenul}
\]
as long as $ \lambda\le\nsize^{1/2} \gmiid(\rri) \le\gmb(\rr)
$.
Similarly for $ (E D_{0}) $ and $ (E D_{1}) $.
\end{pf}

%re7 #&#
\begin{remark}
\label{Rcondpa}
This remark discusses how the presented conditions relate to what is
usually assumed
in statistical literature.
One general remark concerns the choice of the parametric family $
(P_{\thetav}) $.
The point of the classical theory is that the true measure is in this
family, so the
conditions should be as weak as possible.
The viewpoint of this paper is slightly different: whatever family
$ (P_{\thetav}) $ is taken, the true measure is never included,
any model is only an approximation of reality.
From the other side, the choice of the parametric model $ (P_{\thetav
}) $ is
always done by a statistician.
Sometimes some special stylized features of the model
force to include an irregularity in this family.
Otherwise any smoothness condition on the density
$ \ell(y,\thetav) $ can be secured by a proper choice of the family
$ (P_{\thetav}) $.

The presented list also includes the exponential moment conditions $
(ed_{0}) $ and
$ (ed_{1}) $ on the gradient $ \nabla\ell(Y_{1},\thetav) $.
We need exponential moments for establishing some nonasymptotic risk bounds;
the classical concentration bounds require even stronger conditions
that the
considered random variables are bounded.

The identifiability condition $ (\ell\rri) $ is very easy to check
in the usual
asymptotic setup.
Indeed, if the parameter set $ \Theta$ is compact, the Kullback--Leibler
divergence $ \kullbi(\thetav,\thetav^{*}) $ is continuous and
positive for all
$ \thetav\ne\thetav^{*}$, then $ (\ell\rri) $ is fulfilled
automatically with
a universal constant $ \gmi$.
If $ \Theta$ is not compact, the condition is still fulfilled but
the function
$ \gmi(\rri) $ may depend on $ \rri$.
\end{remark}

Below we specify the general results of Sections~\ref{Chglocal}
and~\ref{Chgparam} to the i.i.d. setup.

%s5.1.1 #&#
\subsubsection{A large deviation bound}
This section presents some sufficient conditions ensuring a
small deviation probability for the event
$ \{ \tilde{\thetav} \notin\Thetasi(\rris) \} $ for a fixed~$\rris$.
Below $ \entrlb= \cdimb\dimp$.
We only discuss the case $ \gmi(\rri) \equiv\gmi$.
The general case only requires more complicated notation.
The next result follows from Theorem~\ref{CThittingglrc} with the
obvious changes.

%th5.2 #&#
\begin{theorem}
\label{TLDiid}
Suppose $ (e\rri) $ and $ (\ell\rri) $ with $ \gmi(\rri)
\equiv\gmi$.
If, for $ \rris> 0 $,
%
%e5.12 #&#
%e5.13 #&#
\begin{eqnarray}
\label{cgmi1rriid} \nsize^{1/2} \rris\gmi & \ge& 6 \nunu\sqrt{\xx+
\entrlb} ,
\nonumber
\\[-9pt]
\\[-9pt]
\nonumber
1 + \sqrt{\xx+ \entrlb} & \le& 3 \gmi^{-1} \nunu^{2} \gmiid(
\rris) \nsize^{1/2} , %\label{cgmi2rriid}
\end{eqnarray}
then
%
%e5.14 #&#
\[
%[c]
\P \bigl( \tilde{\thetav} \notin\Thetasi(\rris) \bigr) = \P \bigl( \bigl\|
\mathbf{v}_{0}\bigl(\tilde{\thetav} - \thetav^{*}\bigr)\bigr \| >
\rris \bigr) \le \ex^{-\xx} . \label{PCAxxgliid}\vspace*{-2pt}
\]
\end{theorem}

%re8 #&#
\begin{remark}
\label{RLDiid}
The presented result helps to qualify two important values $ \rris$
and $ \nsize$
providing a sensible deviation probability bound.
For simplicity suppose that $ \gmiid(\rri) \equiv\gmiid> 0 $.
Then the condition \eqref{cgmi1rriid} can be written as
$ \nsize\rris^{2} \gg\xx+ \entrlb$.
In other words, the result of the theorem claims a large deviation
bound for the
vicinity $ \Thetasi(\rris) $ with $ \rris^{2} $ of order $
\dimp/\nsize$.
In classical asymptotic statistics this result is usually referred to as
\emph{root-n consistency}.
Our approach yields this result in a very strong form and for finite samples.\vspace*{-2pt}
\end{remark}

%s5.1.2 #&#
\subsubsection{Local inference}
Now we restate the general local bounds of Section~\ref{SLBuLI} for
the i.i.d. case.
First we describe the approximating linear models.
The matrices $ \mathbf{v}_{0}^{2} $ and $ \IFc$ from conditions
$ (ed_{0}) $, $ (ed_{1}) $ and $ (\ell_{0}) $ determine their
drift and variance
components.
Define
%
%e5.15 #&#
\[
%[c]
\IFb \eqdef \IFc(1 - \rddelta) - \rdomega\mathbf{v}_{0}^{2}
. \label{IFbex}
\]
If $ \tau_{\rd} \eqdef\rddelta+ \fis^{2} \rdomega< 1 $, then
%
%e5.16 #&#
\[
%[c]
\IFb \ge (1 - \tau_{\rd}) \IFc > 0 . \label{IFbIfctau}
\]
Further, $ \DPb^{2} = \nsize\IFb$ and
%
%e5.17 #&#
\[
%[c]
\label{xivbiid} \xivb \eqdef \DPb^{-1} \nabla\zeta\bigl(
\thetav^{*}\bigr) = ( \nsize\IFb )^{-1/2} \sum\nabla
\ell\bigl(Y_{i},\thetav^{*}\bigr) .
\]
The upper bracketing process reads as
%
%e5.18 #&#
\[
%[c]
\Lab\bigl(\thetav,\thetav^{*}\bigr) = \bigl(\thetav-
\thetav^{*}\bigr)^{\T} \DPb\xivb- \bigl\| \DPb\bigl(\thetav-
\thetav^{*}\bigr) \bigr\|^{2}/2 . \label{Labiid}
\]
This expression can be viewed as log-likelihood
for the linear model $ \xivb= \DPb\thetav+ \varepsilonv$
for a standard normal error $ \varepsilonv$.
The (quasi) MLE $ \tilde{\thetav}_{\rd} $ for this model is of the form
$ \tilde{\thetav}_{\rd} = \DPb^{-1} \xivb$.\vspace*{-2pt}

% \alp_{\rd}
% & \eqdef&
% \| \Id_{\dimp} - \IFb^{1/2} \IFm^{-1} \IFb^{1/2} \|_{\infty}
% =
% \lambda_{\max}\bigl( \Id_{\dimp} - \IFb^{1/2} \IFm^{-1} \IFb^{1/2}

%th5.3 #&#
\begin{theorem}
\label{Tconflociid}
Suppose $ (ed_{0}) $.
Given $ \rris$, assume $ (ed_{1}) $, $ (\ell_{0}) $ and $
(\iota) $
on $ \Thetasi(\rris) $, and let $ \rdomega= 3 \nunu  \rhor^{*}
\rris$,
$ \rddelta= \rddelta^{*}\rris$,
and $ \tau_{\rd} \eqdef\rddelta+ \fis^{2} \rdomega< 1 $.
Then the results of Theorem~\ref{TapproxLL} and all its corollaries
apply to the case of
i.i.d. modeling with $ \rups^{2} = \nsize\rris^{2} $. In particular,
on the random set $ \LCS_{\rd}(\rups)
=  \{ \tilde{\thetav} \in\Thetasi(\rris), \| \xivm\| \le
\rups \} $,
it holds
%
%e5.19 #&#
%e5.20 #&#
\begin{eqnarray*}
\| \xivm\|^{2}/2 - \errm(\rups) &\le& L\bigl(\tilde{\thetav},
\thetav^{*}\bigr)  \le \| \xivb\|^{2}/2 + \errb(\rups) ,
\label{vcctttnL}
\\[-2pt]
\bigl\| \sqrt{\nsize\IFb} \bigl( \tilde{\thetav} - \thetav^{*} \bigr) -
\xivb%\DPc^{-1} \nabla L(\thetavs)
\bigr\|^{2} & \le& 2 \spread_{\rd}(\rups) .\vadjust{\goodbreak}
\end{eqnarray*}
The random quantities $ \errb(\rups) $, $ \errm(\rups) $ and $
\spread_{\rd}(\rups) $ follow
the probability bounds of Theorems~\ref{Tdefbound} and~\ref{LLbrevelocm}.
\end{theorem}

Now we briefly discuss the implications of Theorem~\ref{TLDiid} and
\ref{Tconflociid} to the classical asymptotic setup with $ \nsize\to
\infty$.
We fix $ \rris^{2} = C \dimp/\nsize$ for a constant $ C $ ensuring
the deviation bound of Theorem~\ref{TLDiid}.
Then $ \rddelta$ is of order $ \rris$ and the same for $
\rdomega$.
For a sufficiently large $ \nsize$, both quantities are small and, thus,
the spread $ \spread_{\rd}(\rups) $ is small as well; see
Section~\ref{Sdeficincy}.

Further, under $ (ed_{0}) $ condition, the normalized score
%
%e5.21 #&#
\[
%[c]
\label{xiviid} \xiv \eqdef ( \nsize\IFc )^{-1/2} \sum
\nabla\ell\bigl(Y_{i},\thetav^{*}\bigr)
\]
is zero mean asymptotically normal by the central limit theorem.
Moreover, if $ \IFc= \mathbf{v}_{0}^{2} $, then $ \xiv$ is
asymptotically standard
normal.
The same holds for $ \xivb$.
This immediately yields all classical asymptotic results like Wilks
theorem or the Fisher
expansion for MLE in the i.i.d. setup as well as the asymptotic efficiency
of the MLE.
Moreover, our results' bounds yield the asymptotic result for the case
when the parameter
dimension $ \dimp= \dimp_{n} $ grows linearly with $ \nsize$.
Below $ u_{n} = o_{n}(\dimp_{n}) $ means that $ u_{n}/\dimp_{n}
\to0 $ as $
n \to\infty$.
%
%th5.4 #&#
\begin{theorem}
\label{Txiviid}
Let $ Y_{1},\ldots,Y_{\nsize} $ be i.i.d. $ \P_{\thetav^{*}} $ and
let $ (ed_{0}) $, $ (ed_{1}) $, $ (\ell_{0}) $, $ (\iota) $,
$ (e\rri) $ and $ (\ell\rri) $ with $ \gmi(\rri) \equiv\gmi
$ hold.
If $ \nsize> C \dimp_{\nsize} $ for a fixed constant $ C $
depending on constants in the
above conditions only, then
%
%e5.22 #&#
\[
%[c]
\bigl\| \sqrt{\nsize\IFc} \bigl( \tilde{\thetav} - \thetav^{*}
\bigr) - \xiv\bigr\|^{2} = o_{\nsize}(\dimp_{\nsize}),\qquad  2 L\bigl(
\tilde{\thetav},\thetav^{*}\bigr) - \| \xiv\|^{2} =
o_{\nsize
}(\dimp_{\nsize}) . \label{onpasiid}
\]
\end{theorem}
This result particularly yields that
$ \sqrt{\nsize\IFc}  ( \tilde{\thetav} - \thetav^{*} ) $
is nearly
standard normal and $ 2 L(\tilde{\thetav},\thetav^{*}) $ is nearly
$ \chi^{2}_{\dimp} $.

%s5.2 #&#
\subsection{Generalized linear modeling}
\label{SexpexGLM}
Now we consider a generalized linear modeling (GLM) which is often used
for describing
some categorical data.
Let $ \cc{P} = (P_{w}, w \in\Ups) $ be an exponential family with
a canonical
parametrization; see, for example, \citet{mccu1989}.
The corresponding log-density can be represented as
$ \ell(y,w) = y w - d(w) $ for a convex function $ d(w) $.
The popular examples are given by the binomial (binary response,
logistic) model with
$ d(w) = \log ( e^{w} + 1  ) $,
the Poisson model with $ d(w) = e^{w} $ and the exponential model with
$ d(w) = - \log(w) $.
Note that linear Gaussian regression
%considered in Section~\ref{Slinregr}
is a special case with $ d(w) = w^{2}/2 $.

A GLM specification means that every observation $ Y_{i} $ has a
distribution from the
family $ \cc{P} $ with the parameter $ w_{i} $
which linearly depends on the regressor $ \Psi_{i} \in\R^{\dimp} $:
%
%e5.23 #&#
\begin{equation}
%[c]
Y_{i} \sim P_{\Psi_{i}^{\T} \thetav^{*}} . \label{GLLmod}
\end{equation}
The corresponding log-density of a GLM reads as
%
%e5.24 #&#
\[
%[c]
L(\thetav) % =
% \sum\{ Y_{i} w_{i} - d(w_{i}) \}
= \sum \bigl\{
Y_{i} \Psi_{i}^{\T} \thetav- d\bigl(
\Psi_{i}^{\T} \thetav \bigr) \bigr\} . \label{LtvGLM}\vadjust{\goodbreak}
\]
%
%The model \eqref{GLLmod} implies the identity
%
Under $ \P_{\thetav^{*}} $ each observation $ Y_{i} $ follows \eqref
{GLLmod}, in
particular, $ \E Y_{i} = d'(\Psi_{i}^{\T} \thetav^{*}) $.
However, similarly to the previous sections, it is accepted that
the parametric model~\eqref{GLLmod} is misspecified.
Response misspecification means that the
vector $ \fv\eqdef\E\Yv$ cannot be represented in the form
$ d'(\Psi^{\T} \thetav) $ whatever $ \thetav$ is.
The other sort of misspecification concerns the data distribution.
The model \eqref{GLLmod} assumes that the $ Y_{i} $'s are independent
and the marginal distribution belongs to the given parametric family $
\cc{P} $.
In what follows, we only assume independent data having certain
exponential moments.
The target of estimation $ \thetav^{*}$ is defined by
%
%e5.25 #&#
\[
%[c]
\thetav^{*} \eqdef \mathop{\argmax}_{\thetav} \E L(\thetav).
\label{tslinregr}
\]
The quasi-MLE $ \tilde{\thetav} $ is defined by maximization of $
L(\thetav) $:
%
%e5.26 #&#
\[
%[c]
\tilde{\thetav} = \mathop{\argmax}_{\thetav} L(\thetav) =
\mathop{\argmax}_{\thetav} \sum \bigl\{ Y_{i}
\Psi_{i}^{\T} \thetav- d\bigl(\Psi_{i}^{\T}
\thetav \bigr) \bigr\} . \label{LtthetavGLM}
\]
Convexity of $ d(\cdot) $ implies that $ L(\thetav) $ is a
concave function of
$ \thetav$, so that the
optimization problem has a unique solution and can be effectively solved.
However, a closed form solution is only available for the constant
regression or for the
linear Gaussian regression.
The corresponding target $ \thetav^{*}$ is the maximizer of the
expected log-likelihood:
%
%e5.27 #&#
\[
%[c]
\thetav^{*} = \mathop{\argmax}_{\thetav} \E L(\thetav) =
\mathop{\argmax}_{\thetav} \sum \bigl\{ \fs_{i}
\Psi_{i}^{\T} \thetav- d\bigl(\Psi_{i}^{\T}
\thetav\bigr) \bigr\} \label{LEtthetavGLM}
\]
with $ \fs_{i} = \E Y_{i} $.
The function $ \E L(\thetav) $ is concave as well and
the vector $ \thetav^{*}$ is also well defined.

Define the individual errors (residuals)
$ \varepsilon_{i} = Y_{i} - \E Y_{i} $.
Below we assume that these errors fulfill some exponential moment conditions.

\begin{longlist}[$(e_{1})$]
\item[$(e_{1})$]
\textit{There exist some constants $ \nunu$ and $ \gmiid> 0 $,
and for every $ i $
a constant $ \expzeta_{i} $ such that
$ \E ( \varepsilon_{i}/\expzeta_{i}  )^{2} \le1 $ and}
% for all \( |\lambda| \le\gmiid\)
%
%e5.28 #&#
\begin{equation}
%[c]
\log\E\exp ( {\lambda\varepsilon_{i}}/{
\expzeta_{i}} ) \le \nunu^{2} \lambda^{2} / 2, \qquad |
\lambda| \le\gmiid. \label{expzetanunu}
\end{equation}
\end{longlist}

A natural candidate for $ \expzeta_{i} $ is $ \sigma_{i} $ where
$ \sigma_{i}^{2} = \E\varepsilon_{i}^{2} $ is the variance of $
\varepsilon_{i} $;
see Lemma~2.13 of the supplement [\citet{supp}]. %~\ref{LMlambda}.
Under \eqref{expzetanunu}, introduce a $ \dimp\times\dimp$ matrix
$ \VPc$
defined by
%
%e5.29 #&#
\begin{equation}
%[c]
\VPc^{2} \eqdef % \sigma^{-2}
\sum
\expzeta_{i}^{2} \Psi_{i} \Psi_{i}^{\T}
. \label{VPlinregr}
\end{equation}
Condition $ (e_{1}) $ effectively means that each error term
$ \varepsilon_{i} = Y_{i} - \E Y_{i} $ has some bounded exponential moments:
for $ |\lambda| \le\gmiid$, it holds
$ f(\lambda) \eqdef
\log\E\exp ( {\lambda\varepsilon_{i}}/\break{\expzeta_{i}}  )
< \infty$.
This implies the quadratic upper bound for the function $ f(\lambda)
$
for $ |\lambda| \le\gmiid$; see Lemma~2.13 of the supplement [\citet{supp}]. %
In words, condition $ (e_{1}) $ requires a light (exponentially
decreasing) tail for
the marginal distribution of each~$ \varepsilon_{i} $.\vadjust{\goodbreak}

Define also
%
%e5.30 #&#
\begin{equation}
N^{-1/2} \eqdef \max_{i} \sup_{\gammav\in\R^{p}}
\frac{\expzeta_{i} |\Psi_{i}^{\T} \gammav|}{\| \VPc\gammav\|} . \label{CPsiexp}\vspace*{-2pt}
\end{equation}

%le5.5 #&#
\begin{lemma}
\label{LED0EDGLM}
Assume $ (e_{1}) $ and let $ \VPc$ be defined by \eqref{VPlinregr}
and $ N $ by \eqref{CPsiexp}.
Then conditions $ (E D_{0}) $ and $ (E\rr) $ follow from $
(e_{1}) $ with
the matrix $ \VPc$ due to \eqref{VPlinregr} and $ \gmb= \gmiid
N^{1/2} $.
Moreover, the stochastic component $ \zeta(\thetav) $ is linear in
$ \thetav$ and the condition $ (E D_{1}) $ is fulfilled with
$ \rhor(\rr) \equiv0 $.\vspace*{-2pt}
\end{lemma}

\begin{pf}
The gradient of the stochastic component $ \zeta(\thetav) $ of
$ L(\thetav) $ does not depend on $ \thetav$, namely,
$
\nabla\zeta(\thetav)
=
\sum\Psi_{i} \varepsilon_{i}
$
with $ \varepsilon_{i} = Y_{i} - \E Y_{i} $.
Now, for any unit vector $ \gammav\in\R^{p} $ and $ \lambda\le
\gmb$,
independence of the $ \varepsilon_{i} $'s implies that
%
%e5.31 #&#
\[
%[c]
\log\E\exp \biggl\{ \frac{\lambda}{\| \VPc\gammav\|} \gammav^{\T} \sum
\Psi_{i} \varepsilon_{i} \biggr\} = \sum
\log\E\exp \biggl\{ \frac{\lambda\expzeta_{i} \Psi_{i}^{\T} \gammav}{\| \VPc\gammav
\|} \varepsilon_{i} /
\expzeta_{i} \biggr\} . \label{logexplinregr}
\]
By definition,
$ \expzeta_{i} |\Psi_{i}^{\T} \gammav| / \| \VPc\gammav\| \le
N^{-1/2} $ and,
therefore,
$ \lambda\expzeta_{i} |\Psi_{i}^{\T} \gammav| / \| \VPc\gammav
\| \le\gmiid$.
Hence, \eqref{expzetanunu} implies
%
%e5.32 #&#
\begin{equation}
%[c]
\log\E\exp \biggl\{ \frac{\lambda}{\| \VPc\gammav\|} \gammav^{\T} \sum
\Psi_{i} \varepsilon_{i} \biggr\} \le
\frac{\nunu^{2} \lambda^{2}}{2\| \VPc\gammav\|^{2}} \sum\expzeta_{i}^{2}\bigl |
\Psi_{i}^{\T} \gammav\bigr|^{2} = \frac{\nunu^{2} \lambda^{2}}{2} ,
\label{logexpLRED}
\end{equation}
and $ (E D_{0}) $ follows.\vspace*{-2pt}
%
%Therefore, condition \( (E D_{1}) \) is fulfilled automatically
\end{pf}

It only remains to bound the quality of quadratic approximation for the
mean of the
process $ L(\thetav,\thetav^{*}) $ in a vicinity of $ \thetav^{*}$.
An interesting feature of the GLM is that the effect of model misspecification
disappears in the expectation of $ L(\thetav,\thetav^{*}) $.\vspace*{-2pt}

%le5.6 #&#
\begin{lemma}
\label{LttGLM}
It holds
%
%e5.33 #&#
%e5.34 #&#
\begin{eqnarray}
\label{LthetavGLM} - \E L\bigl(\thetav,\thetav^{*}\bigr) % =
% \sum\{ Y_{i} w_{i} - d(w_{i}) \}
&=& \sum \bigl\{ d\bigl(\Psi_{i}^{\T}
\thetav\bigr) - d\bigl(\Psi_{i}^{\T} \thetav^{*}
\bigr) - d'\bigl(\Psi_{i}^{\T}
\thetav^{*}\bigr) \Psi_{i}^{\T} \bigl(\thetav-
\thetav^{*}\bigr) \bigr\}
\nonumber
\\[-8pt]
\\[-8pt]
\nonumber
&=& \kullb ( \P_{\thetav^{*}},\P_{\thetav} ),
\end{eqnarray}
where $ \kullb ( \P_{\thetav^{*}},\P_{\thetav}  ) $ is
the Kullback--Leibler
divergence between measures $ \P_{\thetav^{*}} $ and~$ \P_{\thetav}
$.
Moreover,
%
%e5.35 #&#
\begin{equation}
%[c]
- \E L\bigl(\thetav,\thetav^{*}\bigr) % =
% \frac{1}{2} \sum
% d''(\Psi_{i}^{\T} \thetavd) |\Psi_{i}^{\T} (\thetav- \thetavs)|^{2}
= \bigl\| \DP\bigl(\thetav^{\circ}\bigr) \bigl(\thetav-
\thetav^{*}\bigr) \bigr\|^{2} / 2, \label{ELthetavGLM}
\end{equation}
where $ \thetav^{\circ}\in[\thetav^{*},\thetav] $ and
%
%e5.36 #&#
\[
%[c]
\DP^{2}\bigl(\thetav^{\circ}\bigr) = \sum
d''\bigl(\Psi_{i}^{\T}
\thetav^{\circ}\bigr) \Psi_{i} \Psi_{i}^{\T}
. \label{DGLM}\vspace*{-2pt}
\]
\end{lemma}

\begin{pf}
The definition implies
%
%e5.37 #&#
\[
%[c]
\E L\bigl(\thetav,\thetav^{*}\bigr) % =
% \sum\{ Y_{i} w_{i} - d(w_{i}) \}
= \sum \bigl\{ \fs_{i} \Psi_{i}^{\T}
\bigl(\thetav- \thetav^{*}\bigr) - d\bigl(\Psi_{i}^{\T}
\thetav\bigr) + d\bigl(\Psi_{i}^{\T} \thetav^{*}
\bigr) \bigr\} . \label{LthetavGLM2}
\]
As $ \thetav^{*}$ is the extreme point of $ \E L(\thetav) $, it holds
$ \nabla\E L(\thetav^{*})
= \sum [ \fs_{i} - d'(\Psi_{i}^{\T} \thetav^{*})  ] \Psi_{i} = 0 $
and \eqref{LthetavGLM} follows.
The Taylor expansion of the second order around $ \thetav^{*}$ yields
the expansion
\eqref{ELthetavGLM}.\vadjust{\goodbreak}
\end{pf}

Define now the matrix $ \DPc$ by
%
%e5.38 #&#
\[
%[c]
\DPc^{2} \eqdef \DP^{2}\bigl(
\thetav^{*}\bigr) = \sum d''\bigl(
\Psi_{i}^{\T} \thetav^{*}\bigr) \Psi_{i}
\Psi_{i}^{\T} . \label{DPcGLM}
\]
Let also $ \VPc$ be defined by \eqref{VPlinregr}.
Note that the matrices $ \DPc$ and $ \VPc$ coincide if the model
$ Y_{i} \sim P_{\Psi_{i}^{\T} \thetav^{*}} $ is correctly specified and
$ \expzeta_{i}^{2} = d''(\Psi_{i}^{\T} \thetav^{*}) $.
The matrix $ \VPc$ describes a local elliptic neighborhood of the
central point
$ \thetav^{*}$ in the form
$ \Thetas(\rr) = \{ \thetav\dvtx  \| \VPc( \thetav- \thetav^{*}) \|
\le
\rr\} $.
If the matrix function $ \DP^{2}(\thetav) $ is continuous in this vicinity
$ \Thetas(\rr) $, then the value $ \rddelta(\rr) $ measuring
the approximation quality of $ - \E L(\thetav,\thetav^{*}) $ by the quadratic
function $ \| \DPc(\thetav- \thetav^{*}) \|^{2}/2 $ is small and the
identifiability
condition $ (\LL_{0}) $ is fulfilled on $ \Thetas(\rr) $.

%le5.7 #&#
\begin{lemma}
\label{LapproxDPcGLM}
Suppose that
%
%e5.39 #&#
\begin{equation}
%[c]
\bigl\| \Id_{\dimp} - \DPc^{-1} \DP^{2}(\thetav)
\DPc^{-1} \bigr\|_{\infty} \le \rddelta(\rr) , \qquad \thetav\in\Thetas(\rr) .
\label{IDPcGLM}
\end{equation}
Then $ (\LL_{0}) $ holds with this $ \rddelta(\rr) $.
Moreover, as the quantities $ \rhor(\rr), \errb(\rr), \errm(\rr)
$ vanish,
one can take $ \rdomega= 0 $, leading to the following
representation for
$ \DPb$ and $ \xivb$:
%
%e5.40 #&#
%e5.41 #&#
\begin{eqnarray*}
\label{DPbGLM} \DPb^{2} &=& (1 - \rddelta) \DPc^{2},\qquad \xivb =
(1 + \rddelta)^{1/2} \xiv,
\\[-2pt]
\DPm^{2} &=& (1 + \rddelta) \DPc^{2},\qquad  \xivm = (1 -
\rddelta)^{1/2} \xiv,
\end{eqnarray*}
with
%
%e5.42 #&#
\[
%[c]
\xiv \eqdef \DPc^{-1} \nabla\zeta = \DPc^{-1} \sum
\Psi_{i} (Y_{i} - \E Y_{i}) .
\label{xivGLM}
\]
\end{lemma}

Linearity of the stochastic component $ \zeta(\thetav) $
in the considered GLM implies the important fact that
the quantities $ \errb(\rr), \errm(\rr) $ in the general
bracketing bound~\eqref{LttbLtt} vanish for any $ \rr$.
Therefore, in the GLM case, the deficiency can be defined as the difference
between upper and lower excess and it can be easily evaluated:
%
%e5.43 #&#
\[
%[c]
\spread(\rr) = \| \xivb\|^{2}/2 - \| \xivm
\|^{2}/2 = \rddelta\| \xiv\|^{2} . \label{defGLM}
\]
Our result assumes some concentration properties of the squared norm
$ \| \xiv\|^{2} $ of the vector $ \xiv$.
These properties can be established by general results of Section~\ref{sec1} of
the complement %\ref{Sprobabquad}
under the regularity condition: for some $ \fis$
%
%e5.44 #&#
\begin{equation}
%[c]
\VPc\le\fis\DPc. \label{VPcDPcGLM}
\end{equation}

Now we are prepared to state the local results for the GLM estimation.

%th5.8 #&#
\begin{theorem}
\label{TLCBlinGLM}
Let $ (e_{1}) $ hold.
% and \( \DPc^{2} \succeq\fis^{2} \VPc^{2} \) with \( \fis> 0 \).
Then for %\( \rd= (\rddelta,0) \) with
$ \rddelta\ge\rddelta(\rr) $ %and any \( \thetav\in\Thetas(\rr)
% \Lam(\thetav,\thetavs)
% \le
% L(\thetav,\thetavs)
% \le
% \Lab(\thetav,\thetavs).
%Moreover, for
any $ z > 0 $ and $ \zz> 0 $, it holds
%
%e5.45 #&#
%e5.46 #&#
\begin{eqnarray*}
\P \bigl(\bigl \| \DPc \bigl( \tilde{\thetav} - \thetav^{*} \bigr)\bigr \| > z,
\bigl\| \VPc \bigl( \tilde{\thetav} - \thetav^{*} \bigr) \bigr\| \le\rr \bigr) &
\le& \P \bigl\{ \| \xiv\|^{2} > (1 - \rddelta) z^{2} \bigr\},
\label{pnnRGLM}
\\[-2pt]
\P \bigl( L\bigl(\tilde{\thetav},\thetav^{*}\bigr) > \zz,\bigl \| \VPc
\bigl( \tilde{\thetav} - \thetav^{*} \bigr)\bigr \| \le\rr \bigr) & \le& \P
\bigl\{ \| \xiv\|^{2}/2 > (1 - \rddelta) \zz \bigr\} .
\label{PlLdrrlinGLM}\vadjust{\goodbreak}
\end{eqnarray*}
Moreover, on the set
$ \LCS_{\rd}(\rr) =
\{ \| \VPc ( \tilde{\thetav} - \thetav^{*} ) \| \le\rr,
\| \xivm\| \le\rr\} $, it holds
%
%e5.47 #&#
\begin{equation}
%[c]
\bigl\| \DPc \bigl( \tilde{\thetav} - \thetav^{*} \bigr) - \xiv
\bigr\|^{2} \le \frac{2\rddelta}{1 - \rddelta^{2}} \| \xiv\|^{2} .
\label{DPcttxiGLM}
\end{equation}
\end{theorem}

If the function $ d(w) $ is quadratic, then the approximation error
$ \rddelta$ vanishes as well and the expansion \eqref{DPcttxiGLM}
becomes equality
which is also fulfilled globally, a localization step in not required.
However, if $ d(w) $ is not quadratic, the result applies only
locally and it
has to be accomplished with a large deviation bound.
The GLM structure is helpful in the large deviation zone as well.
Indeed, the gradient $ \nabla\zeta(\thetav) $ does not depend on
$ \thetav$
and, hence, the most delicate condition $ (E\rr) $ is fulfilled
automatically with
$ \gmb= \gmiid N^{1/2} $ for all local sets $ \Thetas(\rr) $.
Further, the identifiability condition $ (\cc{L}\rr) $ easily
follows from
Lemma~\ref{LttGLM}: it suffices to bound from below the matrix
$ \DP(\thetav) $ for $ \thetav\in\Thetas(\rr) $:
%
%e5.48 #&#
\[
%[c]
\DP(\thetav) \ge \gmi(\rr) \VPc,\qquad \thetav\in\Thetas(\rr) .
\label{DPttVPGLM}
\]

An interesting question, similarly to the i.i.d. case, is the minimal radius
$ \rups$ of the local vicinity $ \Thetas(\rups) $ ensuring the
desirable concentration property.
Suppose for the moment that the constants $ \gmi(\rr) $ are all the
same for
different $ \rr$: $ \gmi(\rr) \equiv\gmi$.
Under the regularity condition \eqref{VPcDPcGLM}, a sufficient lower
bound for
$ \rups$ can be based on Corollary~\ref{CThittingglr}.
The required condition can be restated as
%
%e5.49 #&#
\[
%[c]
1 + \sqrt{\xx+ \entrlb} \le 3 \nunu^{2} \gm/\gmi,\qquad  6 \nunu
\sqrt{\xx+ \entrlb} \le \rr\gmi. \label{cgmi2rr}
\]
It remains to note that $ \entrlb= \cdimb\dimp$ and $ \gmb=
\gmiid N^{1/2} $.
So, the required conditions are fulfilled for $ \rr^{2} \ge\rups^{2} = C (\xx+ \dimp) $,
where $ C $ only depends on $ \nunu, \gmi$, and $ \gm$.
%

%s5.3 #&#
\subsection{Linear median estimation}
\label{Slinmedian}
This section illustrates how the proposed approach applies to robust estimation
in linear models.
The target of analysis is the linear dependence of the observed data
$ \Yv= (Y_{1},\ldots,Y_{\nsize}) $ on the set of features
$ \Psi_{i} \in\R^{\dimp} $:
%
%e5.50 #&#
\begin{equation}
%[c]
Y_{i} = \Psi_{i}^{\T} \thetav+
\varepsilon_{i} , \label{regrmod}
\end{equation}
where $ \varepsilon_{i} $ means the $ i $th individual error.
As usual, the true data distribution can deviate from the linear model.
In addition, we admit contaminated data which naturally leads to the
idea of robust
estimation.
This section offers a qMLE view on the robust estimation problem.
Our parametric family assumes the linear dependence \eqref{regrmod}
with i.i.d. errors $ \varepsilon_{i} $ which follow the double
exponential (Laplace)
distribution with the density $ (1/2) e^{-|y|} $.
Then the corresponding log-likelihood reads as
%
%e5.51 #&#
\[
%[c]
L(\thetav) = - \frac{1}{2} \sum\bigl|Y_{i} -
\Psi_{i}^{\T} \thetav\bigr| \label{LthetaLAD}
\]
and $ \tilde{\thetav} \eqdef\argmax_{\thetav} L(\thetav) $ is called
the \emph{least absolute deviation} (LAD) estimate.
In the context of linear regression, it is also called the \emph
{linear median} estimate.
The target of estimation $ \thetav^{*}$ is usually defined by the equation
$ \thetav^{*}= \argmax_{\thetav} \E L(\thetav) $.

It is useful to define the residuals
$ \tilde{\varepsilon}_{i} = Y_{i} - \Psi_{i}^{\T} \thetav^{*}$ and
their distributions
%
%e5.52 #&#
\[
%[c]
P_{i}(A) = \P ( \tilde{\varepsilon}_{i} \in A )
= \P \bigl( Y_{i} - \Psi_{i}^{\T}
\thetav^{*}\in A \bigr) \label{PiAmed}
\]
for any Borel set $ A $ on the real line.
%Define also \( \P_{\nsize} = \nsize^{-1} (P_{1} + \ldots+ P_{\nsize})
If $ Y_{i} = \Psi_{i}^{\T} \thetav^{*}+ \varepsilon_{i} $ is the
true model,
%and if all the \( \varepsilon_{i} \)'s are i.i.d.
then $ P_{i} $ coincides with the distribution of
each $ \varepsilon_{i} $.
Below we suppose that each $ P_{i} $ has a positive density $ \dens_{i}(y) $.

Note that the difference $ L(\thetav) - L(\thetav^{*}) $ is bounded by
$ \frac{1}{2}\sum|\Psi_{i}^{\T}(\thetav- \thetav^{*})| $.
% and condition \( (E) \) is fulfilled automatically.
Next we check conditions $ (E D_{0}) $ and $ (E D_{1}) $.
Denote $ \xi_{i}(\thetav)
= \Ind(Y_{i} - \Psi_{i}^{\T} \thetav\le0) - \bern_{i}(\thetav) $
for $ \bern_{i}(\thetav) = \P(Y_{i} - \Psi_{i}^{\T} \thetav\le
0) $.
This is a centered Bernoulli random variable, and it is easy to check that
%
%e5.53 #&#
\begin{equation}
%[c]
\nabla\zeta(\thetav) = - \sum\xi_{i}(\thetav)
\Psi_{i} . \label{nablazetaLAD}
\end{equation}
This expression differs from the similar ones from the linear and
generalized linear
regression because the stochastic terms $ \xi_{i} $ now depend on $
\thetav$.
First we check the global condition $ (E\rr) $.
Fix any $ \gmiid< 1 $.
Then it holds for a Bernoulli r.v. $ Z $ with
$ \P(Z = 1) = \bern$, $ \xi= Z - \bern$, and $ |\lambda| \le
\gmiid$
%
%e5.54 #&#
%e5.55 #&#
\begin{eqnarray}\label{xideltaB}
\log\E\exp(\lambda\xi) &=& \log \bigl[ \bern\exp\bigl\{ \lambda(1 - \bern) \bigr
\} + (1 - \bern) \exp(- \lambda \bern) \bigr]
\nonumber
\\[-8pt]
\\[-8pt]
\nonumber
& \le& \nunu^{2} \bern(1 - \bern) \lambda^{2} / 2 ,
\end{eqnarray}
where $ \nunu\ge1 $ depends on $ \gmiid$ only.
Let now a vector $ \gammav\in\R^{p} $ and $ \rho> 0 $ be such that
$ \rho|\Psi_{i}^{\T} \gammav| \le\gmiid$ for all $ i=1,\ldots
,\nsize$.
Then
%
%e5.56 #&#
%e5.57 #&#
\begin{eqnarray}\label{logxideltaB}
\log\E\exp\bigl\{ \rho\gammav^{\T} \nabla\zeta(\thetav) \bigr\} & \le&
\frac{\nunu^{2} \rho^{2}}{2} \sum\bern_{i}(\thetav) \bigl\{ 1 -
\bern_{i}(\thetav) \bigr\} % \sigma^{2}\Bigl( \bern_{i}\bigl( \Psi_{i}^{\T} (\thetav- \thetavs)
\bigl|
\Psi_{i}^{\T} \gammav\bigr|^{2}
\nonumber
\\[-8pt]
\\[-8pt]
\nonumber
& \le& \frac{\nunu^{2} \rho^{2}}{2} \bigl\| \VP(\thetav) \gammav\bigr\|^{2} ,
\end{eqnarray}
where %\( \sigma^{2}(\bern) = \bern(1 - \bern) \le1/4 \) and
%
%e5.58 #&#
\[
%[c]
\VP^{2}(\thetav) = \sum\bern_{i}(
\thetav) \bigl\{ 1 - \bern_{i}(\thetav) \bigr\} % \sigma^{2}(\bern_{i}\bigl( \Psi_{i}^{\T} (\thetav- \thetavs) \bigr))
\Psi_{i} \Psi_{i}^{\T} . \label{VPthetaLAD}
\]
Denote also
%
%e5.59 #&#
\begin{equation}
% \VPc^{2}
% &=&
% \VP^{2}(\thetavs)
% =
% \sum\bern_{i}(\thetavs) \bigl\{ 1 - \bern_{i}(\thetavs) \bigr\}
% \Psi_{i} \Psi_{i}^{\T} ,
% \\
\VPc^{2} = \frac{1}{4} \sum
\Psi_{i} \Psi_{i}^{\T} . \label{VP214}
\end{equation}
Clearly, $ \VP(\thetav) \le\VPc$ for all $ \thetav$
and condition $ (E\rr) $ is fulfilled with the matrix $ \VPc$
and $ \gmb(\rr) \equiv\gmb= \gmiid N^{1/2} $ for $ N $ defined by
%
%e5.60 #&#
\begin{equation}
%[c]
N^{-1/2} \eqdef \max_{i} \sup_{\gammav\in\R^{\dimp}}
\frac{\Psi_{i}^{\T} \gammav}{2 \| \VPc\gammav\|} ; \label{NLAD}
\end{equation}
cf. \eqref{logexpLRED}.\vadjust{\goodbreak}

Let some $ \rups> 0 $ be fixed.
We will specify this choice later.
Now we check the local conditions within the elliptic vicinity
$ \Thetas(\rups) = \{ \thetav\dvtx  \| \VPc(\thetav- \thetav^{*}) \|
\le
\rups\} $
of the central point $ \thetav^{*}$ for $ \VPc$ from \eqref{VP214}.
Then condition $ (E D_{0}) $ with the matrix $ \VPc$ and
$ \gmb= N^{1/2} \gmiid$ is fulfilled on $ \Thetas(\rups) $ due to
\eqref{logxideltaB}.
Next, in view of \eqref{NLAD}, it holds
$ |\Psi_{i}^{\T} \gammav| \le2 N^{-1/2} \| \VPc\gammav\| $ for any
vector $ \gammav\in\R^{\dimp} $.
By \eqref{nablazetaLAD},
%
%e5.61 #&#
\[
\nabla\zeta(\thetav) - \nabla\zeta\bigl(\thetav^{*}\bigr) = \sum
\Psi_{i} \bigl\{ \xi_{i}(\thetav) -
\xi_{i}\bigl(\thetav^{*}\bigr) \bigr\}. \label{nzetaLAD}
\]
If $ \Psi_{i}^{\T} \thetav\ge\Psi_{i}^{\T} \thetav^{*}$, then
%
%e5.62 #&#
\[
%[c]
\xi_{i}(\thetav) - \xi_{i}\bigl(
\thetav^{*}\bigr) = \Ind\bigl(\Psi_{i}^{\T}
\thetav^{*}\le Y_{i} < \Psi_{i}^{\T}
\thetav\bigr) - \P \bigl( \Psi_{i}^{\T} \thetav^{*}
\le Y_{i} < \Psi_{i}^{\T} \thetav \bigr) .
\label{xiitts}
\]
Similarly, for $ \Psi_{i}^{\T} \thetav< \Psi_{i}^{\T} \thetav^{*}$
%
%e5.63 #&#
\[
%[c]
\xi_{i}(\thetav) - \xi_{i}\bigl(
\thetav^{*}\bigr) = - \Ind\bigl(\Psi_{i}^{\T}
\thetav\le Y_{i} < \Psi_{i}^{\T}
\thetav^{*}\bigr) + \P \bigl( \Psi_{i}^{\T} \thetav
\le Y_{i} < \Psi_{i}^{\T} \thetav^{*}
\bigr) . \label{xiittsm}
\]
Define $ \bern_{i}(\thetav,\thetav^{*}) \eqdef
| \bern_{i}(\thetav) - \bern_{i}(\thetav^{*})  | $.
Now \eqref{xideltaB} yields similarly to \eqref{logxideltaB}
%
%e5.64 #&#
%e5.65 #&#
\begin{eqnarray*}\label{logxideltaED1}
 && \log\E\exp \bigl\{ \rho\gammav^{\T} \bigl\{ \nabla\zeta(
\thetav) - \nabla\zeta\bigl(\thetav^{*}\bigr) \bigr\} \bigr\}
 \\
% & \le&
&&\qquad\le \frac{\nunu^{2} \rho^{2}}{2} \sum
\bern_{i}\bigl(\thetav,\thetav^{*}\bigr) \bigl|
\Psi_{i}^{\T} \gammav\bigr|^{2}
\\
&&\qquad \le 2 \nunu^{2} \rho^{2} \max_{i \le n}
\bern_{i}\bigl(\thetav,\thetav^{*}\bigr) \| \VPc\gammav
\|^{2} \le \rhor(\rr) \nunu^{2} \rho^{2} \| \VPc
\gammav\|^{2} / 2 ,
\end{eqnarray*}
with
%
%e5.66 #&#
\[
%[c]
\rhor(\rr) \eqdef 4 \max_{i \le\nsize} \sup_{\thetav\in\Thetas(\rr)}
\bern_{i}\bigl(\thetav,\thetav^{*}\bigr) . \label{CEDLAD}
\]
If each density function $ p_{i} $ is uniformly bounded by a constant
$ C $, then
%
%e5.67 #&#
\[
%[c]
\bigl|\bern_{i}(\thetav) - \bern_{i}\bigl(
\thetav^{*}\bigr)\bigr| \le C \bigl| \Psi_{i}^{\T} \bigl(
\thetav- \thetav^{*}\bigr)\bigr | \le C N^{-1/2} \bigl\| \VPc\bigl(
\thetav- \thetav^{*}\bigr)\bigr \| \le C N^{-1/2} \rr.
\label{bebeCtts}
\]

Next we check the local identifiability condition.
We use the following technical lemma.
%
%le5.9 #&#
\begin{lemma}
\label{Lidentmed}
It holds for any $ \thetav$
%
%e5.68 #&#
\begin{equation}
%[c]
- \frac{\partial^{2}}{\partial^{2} \thetav} \E L(\thetav) = \DP^{2}(\thetav)
\eqdef \sum p_{i} \bigl( \Psi_{i}^{\T}
\bigl(\thetav- \thetav^{*}\bigr) \bigr) \Psi_{i}
\Psi_{i}^{\T} , \label{d2ELtt}
\end{equation}
where $ \dens_{i}(\cdot) $ is the density of
$ \tilde{\varepsilon}_{i} = Y_{i} - \Psi_{i}^{\T} \thetav^{*}$.
Moreover, there is $ \thetav^{\circ}\in[\thetav,\thetav^{*}] $
such that
%
%e5.69 #&#
%e5.70 #&#
\begin{eqnarray}\label{EYLpsi}
- \E L\bigl(\thetav,\thetav^{*}\bigr) &=& \frac{1}{2} \sum
\bigl|\Psi_{i}^{\T} \bigl(\thetav- \thetav^{*}
\bigr)\bigr|^{2} \dens_{i}\bigl(\Psi_{i}^{\T}
\bigl(\thetav^{\circ}- \thetav^{*}\bigr)\bigr)
\nonumber
\\[-8pt]
\\[-8pt]
\nonumber
&=& \bigl(\thetav- \thetav^{*}\bigr)^{\T} \DP^{2}
\bigl(\thetav^{\circ}\bigr) \bigl(\thetav- \thetav^{*}\bigr) / 2.
\end{eqnarray}
\end{lemma}
\begin{pf}
Obviously
%
%e5.71 #&#
\[
%[c]
\frac{\partial\E L(\thetav)}{\partial\thetav} = \sum \bigl\{ \P\bigl(Y_{i}
\le\Psi_{i}^{\T} \thetav\bigr) - 1/2 \bigr\}
\Psi_{i} . \label{pELtt}
\]
The identity \eqref{d2ELtt} is obtained by one more differentiation.
By definition, $ \thetav^{*}$ is the extreme point of $ \E L(\thetav)
$.
The equality $ \nabla\E L(\thetav^{*}) = 0 $ yields
%
%e5.72 #&#
\[
%[c]
\sum \bigl\{ \P\bigl(Y_{i} \le
\Psi_{i}^{\T} \thetav^{*}\bigr) - 1/2 \bigr\}
\Psi_{i} = 0 . \label{sumpsivs}
\]
Now \eqref{EYLpsi} follows by the Taylor expansion of the second order at
$ \thetav^{*}$.
\end{pf}
Define
%
%e5.73 #&#
\begin{equation}
%[c]
\DPc^{2} \eqdef \sum\bigl|\Psi_{i}^{\T}
\bigl(\thetav- \thetav^{*}\bigr)\bigr|^{2} \dens_{i}(0)
. \label{DPcLAD}
\end{equation}
Due to this lemma, condition $ (\LL_{0}) $ is fulfilled in $
\Thetas(\rr) $ with
this choice $ \DPc$ for $ \rddelta(\rr) $ from \eqref{IDPcGLM};
see Lemma~\ref{LapproxDPcGLM}.
Moreover, if $ \dens_{i}(0) \ge\fis^{2}/4 $ for $ \fis> 0 $, then
the identifiability condition $ (\AssId) $ is also satisfied.
Now all the local conditions are fulfilled,
yielding the general bracketing bound of Theorem~\ref{TapproxLL} and
all its
corollaries.

It only remains to accomplish them by a large deviation bound, that is,
to specify
the local vicinity $ \Thetas(\rups) $ providing the prescribed
deviation bound.
A sufficient condition for the concentration property is that the expectation
$ \E L(\thetav,\thetav^{*}) $ grows in absolute value with the distance
$ \| \VPc(\thetav- \thetav^{*}) \| $.
We use the representation \eqref{d2ELtt}.
Suppose that for some fixed $ \delta< 1/2 $ and $ \rho> 0 $
%
%e5.74 #&#
\begin{equation}
%[c]
\bigl| \dens_{i}(u)/ \dens_{i}(0) - 1 \bigr| \le\delta,
\qquad |u| \le\rho. \label{piupi0}
\end{equation}
For any $ \thetav$ with $ \| \VPc(\thetav- \thetav^{*}) \| = \rr
\ge\rups$,
and for any $ i = 1,\ldots,n $, it holds
%
%e5.75 #&#
\[
%[c]
\bigl|\Psi_{i}^{\T} \bigl(\thetav-
\thetav^{*}\bigr)\bigr| \le N^{-1/2} \bigl\| \VPc\bigl(\thetav-
\thetav^{*}\bigr)\bigr \| = N^{-1/2} \rr. \label{Psiittsme}
\]
Therefore, for $ \rr\le\rho N^{1/2} $ and any $ \thetav\in
\Thetas(\rr) $
with $ \| \VPc(\thetav- \thetav^{*}) \| = \rr$, it holds
$ \dens_{i} ( \Psi_{i}^{\T} (\thetav^{\circ}- \thetav^{*})
)
\ge(1 - \delta) \dens_{i}(0) $.
Now Lemma~\ref{Lidentmed} implies
%
%e5.76 #&#
\[
%[c]
- \E L\bigl(\thetav,\thetav^{*}\bigr) \ge
\frac{1 - \delta}{2} \bigl\| \DPc\bigl(\thetav- \thetav^{*}\bigr)
\bigr\|^{2} \ge \frac{1 - \delta}{2 \fis^{2}} \bigl\| \VPc\bigl(\thetav-
\thetav^{*}\bigr) \bigr\|^{2} = \frac{1 - \delta}{2 \fis^{2}}
\rr^{2} . \label{pmNrrVPcme}
\]
By Lemma~\ref{Lidentmed} the function $ - \E L(\thetav,\thetav^{*})
$ is
convex. This easily yields
%
%e5.77 #&#
\[
%[c]
- \E L\bigl(\thetav,\thetav^{*}\bigr) \ge
\frac{1 - \delta}{2 \fis^{2}} \rho N^{1/2} \rr \label{fisrr2me}
\]
for all $ \rr\ge\rho N^{1/2} $.
Thus,
%
%e5.78 #&#
\[
%[c]
\rr\gmi(\rr) \ge \cases{ (1 - \delta) \bigl(2 \fis^{2}
\bigr)^{-1} \rr,&\quad $ \mbox{if } \rr\le\rho N^{1/2},$ \vspace*{2pt}
\cr
(1 - \delta) \bigl(2 \fis^{2}\bigr)^{-1} \rho
N^{1/2}, &\quad $\mbox{if } \rr> \rho N^{1/2}.$ } \label{rrgmirrme}
\]
So, the global identifiability condition $ (\LL_{1}) $ is fulfilled if
$ \rups^{2} \ge C_{1} \fis^{2} (\xx+ \entrlb) $ and if
$ \rho^{2} N \ge C_{2} \fis^{2} (\xx+ \entrlb) $
for some fixed constants $ C_{1} $ and $ C_{2} $.

Putting this all together yields the following result.
%
%th5.10 #&#
\begin{theorem}
\label{Tmedregr}
\label{TLCBlinMed}
Let $ Y_{i} $ be independent,
$ \thetav^{*}= \argmax_{\thetav} \E L(\thetav) $, $ \DPc^{2} $ be
given by \eqref{DPcLAD},
and $ \VPc^{2} $ by \eqref{VP214}.
Let also the densities $ \dens_{i}(\cdot) $ of $ Y_{i} - \Psi_{i}^{\T} \thetav^{*}$
be uniformly bounded by a constant $ C $,
fulfill \eqref{piupi0} for some $ \rho> 0 $ and $ \delta> 0 $, and
$ \dens_{i}(0) \ge\fis^{2}/4 $ for all $ i $.
Finally, let $ N \ge C_{2}\rho^{-2} \fis^{2} (\xx+ \dimp) $ for some
fixed $ \xx> 0 $ and $ C_{2} $.
Then on the random set of probability at least $ 1 - \ex^{-\xx} $,
one obtains for $ \xiv\eqdef\DPc^{-1} \nabla L(\thetav^{*}) $ the bounds
%
%e5.79 #&#
\[
%[c]
\bigl\| \sqrt{\DPc} \bigl( \tilde{\thetav} - \thetav^{*} \bigr)
- \xiv \bigr\|^{2} = o(\dimp), \qquad 2 L\bigl(\tilde{\thetav},
\thetav^{*}\bigr) - \| \xiv\|^{2} = o(\dimp) .
\label{onpasmed}
\]
\end{theorem}

\section*{Acknowledgments}
Financial support by the German Research Foundation (DFG) through the
Collaborative
Research Center 649 ``Economic Risk'' is gratefully acknowledged.
Critics and suggestions of two anonymous referees helped a lot in
improving the paper.

%================================ append
%================================

\begin{supplement}%[id=suppA]
\stitle{Some results from the theory of empirical processes\\}
\slink[doi]{10.1214/12-AOS1054SUPP} %[doi,text={...}] - jei reikia
%suskaldyti doi
\sdatatype{.pdf}
\sfilename{aos1054\_supp.pdf}
\sdescription{This part collects some general deviation bounds for
non-Gaussian quadratic forms and for general centered random processes
used in the text.}
\end{supplement}

%
% imsref loaded by akundreckaite, 2012-11-20 12:50:00
% imsref loaded by akundreckaite, 2012-11-20 13:00:58
%

\printaddresses

\end{document}